\documentclass[12pt]{article}
\usepackage{geometry}                % See geometry.pdf to learn the layout options. There are lots.
\usepackage{graphicx}
\usepackage{amssymb}
\usepackage{epstopdf}
\usepackage{amsmath}
\usepackage{amsthm}
\usepackage{amsfonts}
\usepackage{epstopdf}

\newtheorem{theorem}{Theorem}[section]
\newtheorem{lemma}[theorem]{Lemma}
\newtheorem{fact}[theorem]{Fact}
\newtheorem{corollary}[theorem]{Corollary}
\newtheorem{claim}[theorem]{Claim}

\newtheorem{proposition}[theorem]{Proposition}
\theoremstyle{definition}
\newtheorem{definition}[theorem]{Definition}

\theoremstyle{remark}
\newtheorem{remark}[theorem]{Remark}
\numberwithin{equation}{section}

\newcommand{\fip}[1]{{\prec \hspace{-0.25em} #1\hspace{-0.25em}\succ}}
\newcommand{\nip}[1]{\langle #1 \rangle}
\newcommand\RR{{\mathbb R}}

\newcommand{\s}[1]{\langle #1 \rangle}
\newcommand\EE{{\mathbb E}}
\newcommand\PP{{\mathbb P}}

\def\W2{W^{1,2}({\cal O}(M))}

\def\Ai{\mbox{Ai}}

\def\b{\beta}
\newcommand{\m}{\mathfrak{m}} %yyy

\newcommand{\abs}[1]{\left| #1 \right|}

\newcommand{\ra}{\rightarrow}
\newcommand{\bareta}{\bar \eta}
\newcommand\eps{\varepsilon}
\newcommand\one{{\mathbf 1}}
\newcommand{\ws}{w}

\newcommand{\Dws}{\Delta w}

\newcommand{\dws}{\delta w}
\newcommand{\sech}{\mbox{sech}}

\newcommand{\hm}[1]{\mbox{H\"older}(#1)^-}

\newcommand{\stb}{\mbox{\small $\frac{2}{\sqrt{\beta}}$}}
\newcommand{\sfb}{\mbox{\small $\frac{4}{\sqrt{\beta}}$}}
\newcommand{\mm}{\mathfrak m}
\newcommand{\Hb}{{\mathcal H}_{\beta}}
\newcommand{\SAE}{SAE$_\beta$}
\newcommand{\HL}{{\mathcal H}^{L}}
\newcommand{\Hloc}{H_{\rm loc}^1}
\newcommand{\Da}{\triangle\,}

\begin{document}

\title{\Large Beta ensembles, stochastic Airy spectrum, and  a diffusion}
\author{Jos\'e A. Ram\'{\i}rez
%\thanks{Department of Mathematics, Universidad de Costa Rica,
%San Jose  2060, Costa Rica}
\and Brian Rider
%\thanks{Department of Mathematics, University of Colorado at Boulder,
%              UCB 395, Boulder, Colorado 80309}
\and
 B\'alint Vir\'ag
 %\thanks{Departments of Mathematics and Statistics, University of Toronto,
  %     ON, M5S 2E4, Canada
              }

\date{}

\maketitle

\begin{abstract}
We prove that the largest eigenvalues of the beta ensembles
of random matrix theory converge in distribution to the
low-lying eigenvalues of the random Schr\"odinger operator
$-\frac{d^2}{dx^2} + x + \frac{2}{\sqrt{\beta}}
b_x^{\prime}$ restricted to the positive half-line, where $b_x^{\prime}$ is white noise. In doing so we extend the
definition of the Tracy-Widom($\beta$) distributions to all
$\beta>0$, and also analyze their tails.  Last, in a parallel development,
we provide a second
characterization of these laws in terms of a
one-dimensional diffusion. The proofs rely on the associated tridiagonal
matrix models and a universality result showing that the
spectrum of such models  converge to that of their
continuum operator limit.  In particular, we show how Tracy-Widom laws arise
from a functional central limit theorem.

\end{abstract}

\section{Introduction}

For any $\beta > 0$, consider the probability density
function of $\lambda_1\ge  \lambda_2 \ge\dots \ge \lambda_n
\in \RR$ given by
\begin{equation}
\label{betadens}
 \PP_n^{\beta}(\lambda_1, \lambda_2, \dots, \lambda_n) =  \frac{1}{Z_{n}^{\beta}}
e^{ - \beta \sum_{k=1}^n \lambda_k^2/4} \, \prod_{j < k} |
\lambda_j - \lambda_k |^{\beta},
\end{equation}
in which $Z_{n}^{\beta}$ is a normalizing constant. When
$\beta= 1, 2 $ or $4$ this is the joint density of
 eigenvalues for the Gaussian orthogonal,
unitary, or symplectic ensembles, G(O/U/S)E, of random
matrix theory.  For these special values of $\beta$, the above
model is solvable: all finite dimensional correlation functions
may be computed explicitly in terms of Hermite functions, allowing
for startling collection of precise local limit theorems for the random points
(see \cite{Deift} for background).

The law (\ref{betadens}) also describes a
one-dimensional Coulomb gas at inverse temperature $\beta$,
and thus is of physical interest.  In fact, (\ref{betadens}) is intimately
connected to Calogero-Sutherland quantum systems and hence
to Jack polynomials (an important system of multiple orthogonal polynomials).
Still, despite long being the focus of several branches of research,
there are no known forms of the $\beta \neq 1, 2, 4$ correlations
which appear amenable to asymptotics.
The forthcoming  and comprehensive text \cite{Forrester} contains
an excellent account of the previous general-beta developments.

It was therefore welcome news when, based on
\cite{Trotter}, Dumitriu and Edelman \cite{DE1} discovered
the following family of matrix models for all $\beta$.
Let $g_1, g_2, \dots g_n$ be independent Gaussians with
mean $0$ and variance 2.  Let also $\chi_{\beta},
\chi_{2\beta }, \dots, \chi_{(n-1)\beta} $ be independent
$\chi$ random variables indexed by the shape parameter.
Then the $n$ eigenvalues of the tridiagonal matrix ensemble
\begin{equation}
\label{thematrix}
H_n^{\beta} =  \frac{1}{\sqrt{\beta}} \left[ \begin{array}{ccccc}  g_1 & \chi_{(n-1)\beta} &&& \\
\chi_{(n-1)\beta} &  g_2 & \chi_{(n-2)\beta} & &\\
&\ddots & \ddots & \ddots & \\
& & \chi_{2\beta } &  g_{n-1} & \chi_{\beta} \\
& & & \chi_{\beta} &  g_{n}  \\
 \end{array}\right]
\end{equation}
have joint law given by (\ref{betadens}). These are more
specifically referred to as the $\beta$-Hermite ensembles.

We focus on the implications of this discovery to the point
process limits of the spectral edge in the general
$\beta$-ensembles. The distributional limits of the largest
eigenvalues in  $G(O/U/S)E$ comprise some of the most
celebrated results in random matrix theory due to their
surprising importance in physics, combinatorics,
multivariate statistics, engineering, and applied
probability: \cite{BDJ}, \cite{BY} \cite{FS}, \cite{J},
\cite{JS},  and \cite{PS} mark a few highlights. The basic
result is: for  $\beta = 1,2,$ or $4$ and $n \uparrow
\infty$, centered by $2\sqrt{  n}$ and scaled by $
n^{1/6}$, the largest eigenvalue converges
in law to the Tracy-Widom($\beta$) distribution (see
\cite{TW1} and \cite{TW2}), which is given explicitly in
terms of the second Painlev\'e transcendent. There are
allied results for second, third, etc.\ eigenvalues; see
again \cite{TW1} as well as
\cite{Dieng}.

The wide array of models for which Tracy-Widom describes
the limit statistics identifies these laws as important new
probability distributions.  Still, our understanding of
these laws is in its infancy. It remains desirable to
obtain a set of characterizing conditions, like those
classically known for say the Gaussian or Poisson laws. A
description of the limit distribution of the largest
eigenvalues in the general $\beta$-Hermite ensembles is a
first step, providing additional information on the
structure of the three Tracy-Widom laws via the structure
of a one parameter family of distributions in which they
naturally reside.

Towards an edge limit theorem at general $\beta$, Sutton \cite{S} and
 Edelman and Sutton \cite{ES}   present a promising heuristic
argument that the rescaled operators
\begin{equation}
\label{scaledmatrix} \tilde{H}_n^{\beta} =   n^{1/6} \left(2\sqrt{ n}\,I-H_n^{\beta} \right),
\end{equation}
where $I$ is the $n \times n$ identity, should correspond
to
\begin{equation}
\label{H}
   {\mathcal H}_{\b} = - \frac{d^2}{dx^2}  + x  + \frac{2}{\sqrt{\beta}}
   \,b^{\prime}_x
\end{equation}
in the $n \uparrow \infty$, or continuum, limit. Here
$b^{\prime}$ indicates a white noise, and the proposed
scaling of the matrix ensembles follows  the edge scaling
in the known cases. Thus, were it to hold, the above
correspondence would entail that the low-lying eigenvalues
of $\tilde H_n^{\beta}$ converge in law to the those of
$\Hb$. Our first result is a proof of this heuristic.

In Section \ref{s.basic} we give a precise definition of this limiting
``stochastic Airy operator"  (\SAE). For now, let $L^*$ denote the space of functions $f$ satisfying $f(0)=0$
and $\int_0^{\infty} (f^{\prime})^2 + (1+x) f^2
 \,dx < \infty$. Then we say $(\psi,\lambda) \in L^*\times \RR$
is an eigenfunction/eigenvalue pair for $\Hb$ if
$\| \psi \|_2 =1$ and
\begin{equation}
\label{e.eqsdef}
    \psi^{\prime\prime}(x) = \stb \,\psi(x) b^\prime_x +  (x - \lambda)  \psi(x),
\end{equation}
holds in the following integration-by-parts sense,
\begin{equation}
\label{e.intparts} \psi^\prime(x) - \psi^\prime(0) = \stb \psi(x)b_x+
\int_0^x -
  \stb \,b_y \psi'(y)\,dy+\int_0^x(y-\lambda)\psi(y)\,dy,
\end{equation}
where all integrands are products of locally $L^2$
functions. The set of eigenvalues is then a deterministic
function of the random Brownian path $b$.

\begin{theorem} \label{mainthm}
With probability one, for each $k\ge0$ the set of eigenvalues of $\Hb$ has a
well-defined $(k+1)$st lowest element $\Lambda_k$. Moreover,
let $\lambda_{1} \ge \lambda_{2} \ge \cdots$ denote the
eigenvalues of the Hermite $\beta$-ensemble
$H_{n}^{\beta}$. Then the vector
\begin{equation}
\Bigl( {n^{1/6}} ( 2\sqrt{
n}-\lambda_{\b, \ell}  ) \Bigr)_{\ell=1,\dots, k}
\end{equation}
converges in distribution as $n \rightarrow \infty$ to
$(\Lambda_0, \Lambda_1, \cdots, \Lambda_{k-1} )$.
\end{theorem}

We mention that though the $n^{1/6}$ scaling across all beta was
anticipated from the Tracy-Widom results at $\beta = 1, 2,
4$, previous authors had only obtained bounds on this
rate. For instance, \cite{D1} proves that the fluctuations of
$\lambda_{\b, 1} $ are no greater than $n^{1/2}$
while, for even integer values of beta, \cite{DF} contains a
suggestive calculation on the one-point function in a vicinity of the
edge which produces  the correct order.

The proof of Theorem \ref{mainthm} relies on an equivalent  variational formulation of the eigenvalue problem.
As we shall show: with self evident notation,
$$
  TW_{\beta} = \sup_{ f \in L^*, || f||_2 = 1 }  \Bigl\{  \stb \int_0^{\infty} f^2(x) db(x) - \int_0^{\infty} [ f'(x)^2 + x f^2(x)] \, dx \Bigr\} ,
$$
where the stochastic integral is again defined via an integration by parts procedure.  The equality here is in law, and one may view this
as a definition of $TW_{\beta}$, independent of any random matrix theory developments.  This variational
approach actually  allows for a type of universality result, and we provide rather general conditions under which
the spectrum of random tridiagonal models converge to that of their
natural continuum operator limit. In particular, one can obtain the
Tracy-Widom laws by considering tridiagonal matrix
versions of $\eqref{H}$ far simpler than
\eqref{scaledmatrix}.

Our next
theorem gives yet another characterization of the limiting spectrum in terms of the explosion
probability of the one-dimensional diffusion $x \mapsto
p(x)$ defined by the It\^o equation
\begin{equation}
\label{diffusion}
   dp\,(x) = \stb \;d b_x + (x -p^2(x)) \, dx.
\end{equation}

\begin{theorem}
\label{difthm} Let $\kappa(x,\cdot)$ be the distribution of
the first passage time to $-\infty$ of the diffusion $p(x)$
when started from $+\infty$ at time $x$.   Then
\begin{eqnarray*}
\PP(\Lambda_0>\lambda)&=&\kappa(-\lambda,\{\infty\}), \qquad \mbox{ and, for }k\ge 0,\\
\PP(\Lambda_{k}<\lambda)&=&
\int_{\RR^{k+1}}\kappa(-\lambda,dx_1)\kappa(x_1,dx_2)\ldots
\kappa(x_{k},dx_{k+1}).
\end{eqnarray*}
\end{theorem}
Even in the well-understood $\beta  =1,2,$ and $4$ settings, any such simple Markovian description of $TW_{\beta}$ is novel. Further, a variant of \eqref{diffusion} is shown to
describe joint laws in Proposition \ref{p.blowups}.

Theorem \ref{difthm} combined with the variational picture
leads to our final result on the shape of
the general beta laws, so far only known for $\beta=1,2,4$.

\begin{theorem}
\label{shapethm}
 With  $TW_{\beta} = - \Lambda_0(\b)$, for $ a \uparrow \infty$ it holds
\begin{eqnarray*}
          \PP \Bigl(  TW_{\beta} > a \Bigr)  &=&  \exp{ \Bigl(-  \frac{2}{3} \beta  a^{3/2} (1 + o(1) )  \Bigr)}, \mbox{ \ \ \ and}\\
  \PP \Bigl(  TW_{\beta} <  -a \Bigr)  &=&  \exp{ \Bigl(-  \frac{1}{24} \beta  a^{3} (1 + o(1) )
  \Bigr)}.
\end{eqnarray*}
\end{theorem}

As indicated, the above discussion carries over to a large class of
$\beta$-ensembles. For example, let $W$
be an $n \times \kappa$ matrix comprised of independent
standard real, complex, or quaternion Gaussians. Viewing
$W$ as $\kappa$ random vectors,  the sample covariance
matrix $W W^{\dagger}$ $-$ known as the Laguerre, or Wishart ensemble
for $\beta=1,2,4$  $-$ plays an important role in mathematical
statistics.

Johansson, in the complex case \cite{J}, and Johnstone, in the
real case \cite{JS}, showed that the largest eigenvalues
tend to the $\beta = 2$ and $\beta = 1$ Tracy-Widom laws
whenever
 $\kappa/n \rightarrow \vartheta \in (0, \infty)$.  Later, El Karoui \cite{EK} showed the same result holds even
 if $\kappa/n \rightarrow \infty$ or $0$ if both $n,\kappa \rightarrow \infty$, this regime being important in applications.
Our proof handles all regimes simultaneously for all
$\beta$.

To explain, consider now the joint density on points $\lambda_1,\dots, \lambda_n \in \RR^+$,
\begin{eqnarray}\label{e.wishart.pd}
  \PP_{n, \kappa}^{\beta}
  (\lambda_1, \lambda_2, \dots, \lambda_n) =  \frac{1}{Z_{n, \kappa}^{\beta}}   \prod_{j < k} | \lambda_j - \lambda_k|^{\beta}
   \times \prod_{k=1}^n  \lambda_k^{\frac{\beta}{2}(\kappa - n+ 1) - 1}  e^{ - \frac{\beta}{2}  \lambda_k} .
\end{eqnarray}
When $\kappa$ is an integer and $\beta = 1$ or $2$ this
is the joint law of Laguerre eigenvalues just described (granting $\kappa \ge n$ which can be assumed
with no loss of generality by the obvious duality). Notice though that the above law is sensible for
any real $\kappa> n-1$ and  $\beta>0$. As an application
of our universal limit theorem, we also show:

\begin{theorem}
\label{Lagthm} Let $\lambda_{1}\ge \lambda_2\ldots$ denote
the ordered $\beta$-Laguerre ``eigenvalues'' \eqref{e.wishart.pd},
and set
\begin{equation}
\label{Wscales}
   \mu_{n, \kappa} =    ( \sqrt{n} + \sqrt{\kappa} )^{2},  \ \mbox{ and } \
  \sigma_{n, \kappa} =
         \frac{ (\sqrt{n \kappa})^{1/3}}{ (\sqrt{n} + \sqrt{\kappa})^{4/3}
         }.
\end{equation}
Then for any $k$, as $n \rightarrow \infty$ with arbitrary
$\kappa=\kappa_n>n-1$ we have
\[
     \Bigl(     \sigma_{n, \kappa}
 ( \mu_{n, \kappa}- \lambda_{\ell}    )  \Bigr)_{\ell = 1, \dots, k}
  \Rightarrow  \Bigl(\Lambda_0 , \Lambda_1, \dots, \Lambda_{k-1}
  \Bigr).
\]
\end{theorem}

In summary, we have a fairly complete characterization of the general beta random matrix ``soft-edge".
Within the context of log-gas type ensembles, the present is the first paper to rigorously establish local fluctuation
theorems of any type away from the classical $\beta = 1,2$ or $4$ exponents on the Vandermonde interaction
component of the density. Since its appearance, there have been several further advances in the
general beta picture: \cite{KS} and \cite{VV1} provide descriptions of the bulk spectrum, with \cite{RR} considering
the so-called ``hard-edge".  The second and third papers rely in part on the results here, and based on \cite{VV1}
those authors have gone on to prove a refined form of Dyson's conjectures for the bulk eigenvalue spacing
for all $\beta > 0$ \cite{VV2}.

Once again,  Section 2 discusses the definition and basic
properties of \SAE.  The diffusion
connection and Theorem \ref{difthm} are detailed in Section
3. Section 4 establishes, the tail bounds, Theorem \ref{shapethm}.
In Section 5
we prove a general result, Theorem
\ref{weak},  which provides weak
 conditions under which the lowest eigenvalues of tridiagonal random matrices of type discrete Laplacian plus
potential converge to the corresponding eigenvalues of
continuum operator limit.  It is anticipated that a result
of this type will be of future importance in investigations
of universality  in random matrix theory, and in Section 5
it is employed to prove Theorems \ref{mainthm} and
\ref{Lagthm}.

\section{Basic properties of the Stochastic Airy Equation}
\label{s.basic}

\subsection*{Definition of the Stochastic Airy Operator}

We use the usual Schwartz distribution theory. Recall
that the space of distributions $D=D(\RR^+)$ is the
continuous dual of the space $C_0^\infty$ of all smooth
compactly supported test functions under the topology of
uniform-on-compact convergence of all derivatives.
Recall as well that all
continuous functions $f$ and their formal derivatives are
distributions. They act on $C_0^\infty$ via integration by
parts,
$$\fip{\varphi,f^{(k)}}:= (-1)^k\int
f(x)\varphi^{(k)}(x)dx,
$$
where the latter is clearly defined. For instance, $b'$, the
formal derivative of Brownian motion, is a random
distribution, as $b$ is a random continuous function. The
notation $\fip{\cdot,\cdot}$ distinguishes the above from
an $L^2$ inner product $\nip{\cdot,\cdot}$.

Introduce $\Hloc$,  the space of functions $f:\RR^+\to \RR$
for which $f^{\prime} \one_I\in L^2$ for any compact set $I$. {\SAE}
is then  well defined as a random
linear map $\Hloc\to D$, sending $f$ to the distribution
$$
\Hb f=-f^{\prime \prime}+ x f  + \stb f b^{\prime}.
$$
As $D$ is only closed under multiplication by smooth
functions, one must make sense  $f b'$ as an element of
that space. Stieltjes integration by parts prompts $\int_0^y f
b'\,dx: = -\int_0^y bf^{\prime}\,dx + f(y)b_y-f(0)b_0$. The latter is a
continuous function of $y$, and we define $fb^{\prime}$ as its
derivative.

\subsection*{Eigenvalues and eigenfunctions}

We will consider the eigenvalues from two points of view,
\begin{itemize}
\item[(i)] as the solutions of
$\Hb f=\lambda f$ with given ``boundary conditions'', or

\item[(ii)]  as solutions of the usual variational problem.
\end{itemize}

The first approach is intimately tied to the Riccati
transformation, which will be our main tool for analyzing
solutions. The second will be  useful in
obtaining bounds on the eigenvalue distributions for \SAE. Their
equivalence  for \SAE\ is established over the course of the
next two subsections.

Recall the Hilbert space  $L^*$ defined via the norm
$$
\|f\|_*^2=\int_0^{\infty} (f^{\prime})^2 + (1+x) f^2 \,dx,
\qquad L^*=\{f:f(0)=0,\|f\|_*<\infty\}.
$$
The framework just introduced makes it self-evident how to
define the $L^*$ eigenfunctions and eigenvalues of \SAE\ in
the sense of (i).

\begin{definition} The eigenvalues and eigenfunctions of $\Hb$ are the
pairs $(f,\lambda)\in L^* \times \RR$ satisfying
$$\Hb f=\lambda f,$$
where both sides are interpreted as distributions ($L^*\subset C^0(\RR^+)\subset D$, where $C^0(\RR^+)$ are those continuous functions
on the half-line, vanishing at infinity). The $(k+1)$st smallest
point in this set, if it
exists, will be denoted $\Lambda_k$.
\end{definition}

The $f\in L^*$  requirement allows us to avoid
technicalities.
Proposition \ref{p.supexp} shows
that it can be relaxed to $f\in C^0\cap L^2$, or even further.

If we rewrite $ \Hb f=\lambda f$  as
\begin{eqnarray}\label{e.ode}
f''= (  x - \lambda + \stb b' )f,
\end{eqnarray}
the right hand side is a derivative of a continuous function,
and thus $f'$ can be taken to be
continuous (it is defined in an a.e.\ sense). In this way we
arrive at two equivalent formulations of the eigenvalue
problem. One is the coincidence of the integrated versions
of both sides of \eqref{e.ode}. The other is the
equivalence of the two sides as distributions. The first
reproduces the definition  \eqref{e.eqsdef} from the introduction
$$
f'(x) - f'(0) =\int_0^x  [  (y-\lambda) f - \stb  b_y f'(y) ]  \,dy \; +   \stb f(x) b_x,
$$
showing at once that  $f'$ inherits the $\hm{1/2}$
continuity properties of $b$, i.e., that it is
H\"older$(1/2-\eps)$ continuous for all $\eps>0$. Thus
$f\in \hm{3/2}$. The second,  weak definition,
will be useful in the variational analysis.  It reads
\begin{equation}
\int \varphi'' f \,dx = \int(x-\lambda) \varphi f \, dx+
\int \stb
 \left[\int_0^x  b_yf'(y)dy-b_xf(x)\right]\varphi' \,dx,
\end{equation}
and
takes the form
\begin{equation}\label{e.weak}
\int \varphi'' f \,dx =\int(x-\lambda)  f\varphi-
\stb bf'\varphi- \stb bf\varphi'\,dx
\end{equation}
after an integration by parts.

Before proceeding to the variational approach,
we register some simple facts about $L^*$.

\begin{fact}\label{f.lstar} Any $L^*$ bounded sequence has
a subsequence  $f_n$ that converges to some $f \in L^{*}$
in each of the following
ways: (i) $f_n \to_{L^2} f$ (ii) $f_n'\to f'$ weakly in
$L^2$, (iii) $f_n\to f$ uniformly-on-compacts and (iv) $f_n
\to f$ weakly in $L^*$.
\end{fact}
\begin{proof} Convergence modes (ii) and (iv) are simple applications of the
Banach-Alaoglu theorem.  (iii) stems from the familiar estimates
$f_n^2(x) \le 2 \| f_n \|_2 \| f_n^{\prime} \|_2$  and
 $ | f_n(y)-f_n(x) | = | \int
f_n'\one_{[x,y]}\,dz | \le\|f_n'\|_2|x-y|^{1/2}$, showing that the
sequence is uniformly equicontinuous on compacts. Last, (iii)
implies $L^2$ convergence locally, while the bound
$\sup_n \int xf_n^2\,dx < \infty $ produces the uniform integrability required for (i).
\end{proof}

\subsection*{The quadratic form $\fip{f,\Hb f }$}

The quadratic form $\fip{ f,\Hb f
}$ typically associated with self-adjoint operators is
already defined for test functions $f\in C_0^\infty$. As
a path to generalization, we proceed by decomposing the
fluctuation term via
$$
b=\bar b+ (b-\bar b), \qquad \bar b_x =\int_x^{x+1} b_y dy,
\mbox{\ \  and so \ \ }\bar b'_x=b_{x+1}-b_{x}.
$$
The idea is to smooth out the noise and then later control the difference
between the noise and its mollification.   We have: for  any $f\in
C_0^\infty$,
\begin{eqnarray}
\fip{f^2, b'}
 &=&
   \int_0^{\infty} f^2(x)\, \bar b'_x\,dx +    \int_0^{\infty} f'(x) f(x) (\bar b_x-
b_x)\, dx.
\label{e.quadraticform}
\end{eqnarray}
The necessary extension requires this object to be finite
over
$L^*$.
Step one is a simple bound on the Brownian paths.

\begin{lemma}\label{lemma:GrowthOfDiffBM}
For $b_x$, $ x > 0$ a Brownian motion, there is a random
constant $C<\infty$ so that
\begin{equation}
\sup_{x>0} \sup_{0 < y \leq 1} \frac{\abs{b_{x+y} -
b_x}}{\sqrt{\log{( 2+ x)}}} \leq C \quad a.s.
\label{GrowthOfDifBM}
\end{equation}
As a consequence, $|\bar b'(x)|\vee|\bar b_x-b_x|\le C
\sqrt{\log( 2+ x)}$.
\end{lemma}

\begin{proof} By the triangle
inequality $ \abs{b_{x+y} - b_x}\leq \abs{b_{x+y} -
b_{\lfloor x+y\rfloor}} + \abs{b_{\lfloor x+y\rfloor} -
b_{\lfloor x \rfloor}} + \abs{b_{x} - b_{\lfloor x\rfloor
}} $ it is enough to show that
\[
    \sup_{n > 1} \frac{X_n}{\sqrt{\log{n}}} \leq C,\]
where the $X_n = \sup_{0 < y \leq 1} \abs{b_{n+y} - b_n}$
are independent and identically distributed. Further,
$P(X_n
> a) = 2 \sqrt{2/\pi} \int_a^{\infty} e^{-m^2/2} dm$ $ \le 2 a^{-1}
e^{-a^2/2}$, courtesy of the reflection principle. A
Borel-Cantelli argument completes the proof.
\end{proof}

Using the lemma, we may certainly bound $  |\bar b'_x| $ by $ C(1+x)$ and
$|\bar b_x-b_x|$ by $\sqrt{C(1+x)}$.
An application of Cauchy-Schwarz  in \eqref{e.quadraticform} then yields

\begin{proposition}\label{l.lstarbound} There is the bound $|\fip{ f^2 , b'}|\le C\|f\|_*^2$ for $C = C(b) < \infty$ a.s.  In particular, $\fip{
f,\Hb f}$ is defined for all $f\in L^*$;  $\fip{f,\Hb
g}$ can be defined by polarization, and is a continuous
symmetric bilinear form $(L^*)^2\to\RR$.
\end{proposition}

\subsection*{Variational characterization}

We are at last in position to define the ground state energy $\tilde \Lambda_0=\tilde
\Lambda_0(\beta, \omega)$ of $\Hb$ variationally.  As anticipated,
we will show that  $\tilde \Lambda_0=\Lambda_0$, the lowest point
in the set of eigenvalues of $\Hb$. Setting
\begin{equation}
\tilde \Lambda_0:= \inf  \left\{ \fip{f,\Hb f}  \;:\; f \in
L^* \mbox{ with }   f(0) = 0, \|f\|_2=1 \right\},
\label{variationalProbH}
\end{equation}
we begin by proving that $\tilde \Lambda_0>C_2(\omega)>-\infty$
a.s.  From now on, random constants are denoted  by $C_\cdot$, and
deterministic ones by $c_\cdot$.

\begin{lemma}\label{Compactness}
\label{l.mainbound} There are constants $c_1,C_2,C_3$ so that
a.s. for all $f\in L^*$
 \begin{equation}\label{ineq:Compact}
 c_1\|f\|_*^2-C_2\|f\|_2^2\le \fip{ f,\Hb f} \le
 C_3\|f\|_*^2.
 \end{equation}
\end{lemma}

\begin{proof}
The upper bound is from Proposition \ref{l.lstarbound}. For the
lower bound, we repeat the definition
\begin{eqnarray*}
\fip{ f,\Hb f }
 &=& \int_0^{\infty} (f')^2 \, dx
 + \int_0^{\infty} xf^2 dx \hspace{8em}\\
 && \lefteqn{+  \stb \int_0^{\infty} f^2(x)\, \bar b'_x\,dx +  \sfb \int_0^{\infty} f'(x) f(x) (\bar b_x-
b_x)\, dx.}
\end{eqnarray*}
The first line on the right
amounts to  $\|f\|_*^2-\|f\|_2^2$, and  it suffices to show
that the terms $B_1, B_2$ in the second line do not ruin
the picture and satisfy $B_1+B_2\ge
-c_1\|f\|_*^2-C\|f\|_2^2$ with $c_1<1$.

Lemma \ref{lemma:GrowthOfDiffBM} provides the  bounds  $\stb |\bar
b'_x|\le c_2(C+x)$ and $ \stb |\bar b_x-b_x|\le c_2\sqrt{C+x}$
for arbitrarily small $c_2$ and some random $C=C(c_2,b)$.
Thus, $|B_1|\le  c_2\|f\|^2_*+c_2C\|f\|_2^2$, and for the
 second term we have
\begin{eqnarray*}
|B_2|\le \int_0^{\infty} |f'(x) f(x) |c_2\sqrt{C+x}\, dx
\le c_2\|f'(x)\|_2^2+c_2\int_0^\infty (C+x) f^2(x)dx.
\end{eqnarray*}
In particular,  $|B|\le 3c_2\|f\|_*^2+C'\|f\|_2^2$, where again $c_2$
may be taken as small as needed.
\end{proof}

\begin{corollary}
The infimum in the variational problem
\eqref{variationalProbH} is attained at an eigenfunction $f_0$
of $\Hb$ with eigenvalue $ \Lambda_0$.
\end{corollary}

\begin{proof}
Again
by Lemma \ref{lemma:GrowthOfDiffBM}, for any
$\eps>0$ the bounds $|\bar b'_x|\le \eps(1+x)$ and $\bar
|b_x-b_x|\le \eps\sqrt{1+x}$ hold for all $x>X(b,c)$. As in
the proof of Lemma \ref{l.mainbound}, it follows
\begin{equation}\label{e.fHftails}
\fip{f,\Hb f} = \|f\|_*^2-\|f\|_2^2 +  \stb \int_0^{X} f^2(x)\,
\bar b'_x\,dx +  \sfb \int_0^{X} f'(x) f(x) (\bar b_x- b_x)\,
dx+\mathcal{E},
\end{equation}
where the error term satisfies $|{\mathcal E}|\le
\eps\|f\|_*^2$.

Next choose a minimizing sequence $f_n\in L^*$:  $||f_n||_2 = 1$  with
 $\fip{f_n,\Hb f_n} \rightarrow \tilde \Lambda_0$. By
Lemma \ref{l.mainbound}, $\|f_n\|_*<B$ for some  $B=B(b)$,
and by Fact \ref{f.lstar}, we can find a subsequence along
which $f_n\to f_0$ uniformly on compacts, in $L^2$, and also
weakly in $H^1$. Terms 2, 3, 4 on the
right hand side of \eqref{e.fHftails} then converge to their evaluations at $f$, while
term 1 plainly satisfies $\|f\|^2_*\le \liminf \|f_n\|^2_*$.
Letting $\eps\to 0$ yields $\fip {f_0,\Hb
f _0}\le \tilde \Lambda_0$, with the opposite inequality holding by
definition.

To complete the picture,  taking the functional derivative
$\frac{d}{d\eps} \fip {f_0+\eps\varphi,\Hb (f_0+\eps
\varphi)}|_{\eps=0}$ of the variational problem
(\ref{variationalProbH}) in directions $\varphi \in
C_0^\infty$ shows that the selected minimizer satisfies $\Hb
f=\tilde{\Lambda}_0 f_0$ in the sense of distributions.  It is therefore an
eigenfunction and $\tilde \Lambda_0 = \Lambda_0$.
\end{proof}

The above formulation may now be used to define  higher
eigenvalues in the expected manner: $\tilde \Lambda_1$ arising
from restricting the class of potential minimizers to be
perpendicular to $f_0$, and so on.  This defines a
sequence of eigenvalues $\tilde \Lambda_k$, $k\ge 0$, each of which
is finite.  That $ \cdots <  \tilde \Lambda_k <  \tilde \Lambda_{k+1} < \cdots $ is proved
later, see Proposition \ref{ToInfty}.  The standard trick with the  bilinear form defined in Lemma
\ref{l.lstarbound} yields  the $L^2$-orthogonality of the corresponding
eigenfunctions  (and, as a by-product, the  uniqueness of the ground state
found above).  In short:

\begin{lemma}
The $(k+1)$st lowest element $\Lambda_k$ in the set of
eigenvalues of $\Hb$ exists and in fact $\Lambda_k=\tilde
\Lambda_k$.
\end{lemma}

\begin{remark}
Lemma \ref{lemma:GrowthOfDiffBM}, as used in Lemma
\ref{l.mainbound}, is related to the second condition of
our general convergence result, Theorem \ref{weak}, below.
In addition, a discrete version of this estimate is employed in the
context of the $\beta$-ensembles converging to \SAE.
\end{remark}

\section{Riccati transform and diffusion formulas}
\label{s.Ricatti}

The {\em Riccati map} is a classical tool in the study
one-dimensional random Schr\"odinger spectra.  Its use
dates back to Halperin \cite{Hlp} who computed the density
of states for $-d^2/dx^2 + b^{\prime}_x$, though see also
\cite{FM}.  For other applications to local statistics
such as the ground state energy, see \cite{CM},
\cite{CRR}, and \cite{M1}.

Return to the eigenvalue problem \eqref{e.eqsdef}
\begin{equation}\label{e.eqsdefrepeat}
    \psi^{\prime\prime}(x) = \stb \,\psi(x) b^\prime_x +  (x - \lambda)  \psi(x),
\end{equation}
understood in  the integration-by parts sense. The Riccati transform is
simply the logarithmic derivative $ p(x) = \psi^{\prime}(x)
/ \psi(x)$. This turns \eqref{e.eqsdefrepeat} into a first
order differential equation: $p(0)=\infty$, and
\begin{equation}
\label{Riccati}
    p'(x)=  x-\lambda-p^2(x)-\stb b'(x),
\end{equation}
understood in the same way.

Solutions of (\ref{Riccati})  may blow up (to $-\infty$) at finite times,
as will happen whenever $\psi$ vanishes. In this case $p$ is
immediately restarted at $+\infty$ at that time-point in order
to continue the
solution corresponding to \eqref{e.eqsdefrepeat}.  It is convenient
to think of $p$ as taking values in the disjoint union of countable copies of the
reals, $\RR_0$,$\RR_{-1}$,$\RR_{-2} \ldots$ Points $(n,x)$
in this space are ordered lexicographically, though we
sometimes refer to these points by their second coordinate
$x$.  A natural topology on this space is provided by the
two-point compactification of each copy of the reals, glued
together at the endpoints so as to respect the
lexicographic ordering. (This ordering and topology can
also be defined by considering the evolution of
$\arg(\psi'+i\psi)$ as a continuous, real-valued function,
and applying the tangent map.)

It is not hard to verify the following.

\begin{fact}
The solution $p_\lambda(x)=p(x,\lambda)$ of
\eqref{Riccati} is unique and increasing in $\lambda$ for
each $x$. It is also decreasing in $x$ at each blowup (or ``explosion").
Moreover, the function $p$ is continuous when the image
space is considered in the topology discussed above.
\end{fact}

Next consider the truncation $\HL$ of $\Hb$, defined on
the finite interval $[0,L]$ with Dirichlet ($\psi=0$)
boundary conditions at both endpoints.

\begin{lemma}\label{l.finiteRiccati}
Fix $\lambda$, and denote  $(-n,y)=p(L,\lambda)$. Then the number
$n$ of blowups of $p(x,\lambda)$ to $-\infty$  on $[0,L]$
equals the number of eigenvalues of $\HL$ at most $\lambda$.
\end{lemma}

\begin{proof}
We provide just a sketch.
First,  $\lambda$ is an eigenvalue of $\HL$ if and only
if $p_\lambda$ blows up to $-\infty$ at the endpoint $L$.
For large negative  $\lambda$, there are no blowups for
any given noise path. As
 $\lambda$ increases, continuity and monotonicity implies that
existing blowups move towards the beginning of the interval $-$
new ones can only appear at the endpoint $L$. At those
$\lambda$ we have a new eigenvalue and the claim follows.
\end{proof}

To extend the picture to the full line we need the following.

\begin{lemma} As $L\to\infty$ the first $k$ eigenvalues
$\Lambda_{L,0},\ldots \Lambda_{L,k-1}$ of $\HL$ converge to the first $k$
eigenvalues of $\Hb $.
\end{lemma}

\begin{proof}
A trivial modification of the proof of Lemma \ref{l.liminf}
together with Fact \ref{f.lstar} shows that $\liminf
\Lambda_{L,k}\ge\Lambda_k$. Next, for an inductive proof,
assume that the $\Lambda_{L,{\ell}}\to \Lambda_\ell$ for
$\ell<k$.

Let $f_k^\eps$ be a function of compact support
$\eps$-close to $f_k$ in $L^*$. Let also
$$g_L=g_{L,k}=f_k^\eps-\sum_{\ell=0}^k \nip{
f_k^\eps,f_{L,\ell}}f_{L,\ell}.$$
The hypothesis entails that $f_{L,\ell}\to_{L^2}f_\ell$. So, for large $L$, each
coefficient in the above sum is bounded by $2\eps$, and
$g_L$ will  be $c\eps$-close to $f_k$ in $L^*$. Then, by the
variational characterization we have
$$
\limsup_{L\to\infty} \Lambda_{L,k} \le \limsup_{L\to\infty}
\frac{\fip{ g_L,\Hb g_L}}{\langle g_L,g_L\rangle },
$$
since $f_k^{\eps}$ is eventually supported on $[0,L]$. Then,
as $\eps\to 0$, the right hand side converges to $\fip{
f_k,\Hb f_k}/{\langle f_k,f_k\rangle }=\Lambda_k$.
\end{proof}

Taking the $L\to\infty$ limit of the claims of Lemma
\ref{l.finiteRiccati} then yields

\begin{proposition}\label{p.blowups} Let $N(\lambda)$ be
the number of blowups of the equation \eqref{Riccati} to
$-\infty$. Then for almost all $\lambda$, $N(\lambda)$
equals the number of eigenvalues of $\Hb$ at most
$\lambda$. In other words,  the cadlag version of $N(\lambda)$ is $\Hb$'s
eigenvalue counting function.
\end{proposition}

Of course, for any fixed $\lambda$, the Riccati equation may be
taken in the It\^o sense,
 \begin{equation}\label{e.Riccatidiffusion}
   dp(x) =  -\stb \, d b_x + (x - \lambda - p^2(x))\, dx,
 \end{equation}
which is to say that $p=p_\lambda = (\log \psi)^{\prime}$
performs the indicated diffusion, restarted at
$+\infty$ instantaneously after each explosion to $-\infty$.
The content of the above is that the total explosion
count equals the count  of eigenvalues $\le \lambda$.

\begin{proof}[Proof of Theorem \ref{difthm}]
The strong Markov property of the motion (\ref{e.Riccatidiffusion})
implies that the sequence of explosion times, $\mm_0=0,\mm_1,\mm_2, \ldots$ is
itself a Markov process. Let $\kappa(x,\cdot)$ be the
distribution of the first such time of
$p_0(x)$ under $\PP_{(\infty,x)}$, that is, when started from $\infty$ at
time $x$. This law is supported on
$(x,\infty]$ with a point mass at $\infty$.
By the preceding,
\begin{eqnarray*}
\PP(\Lambda_{k-1}<\lambda)&=&\PP_{(\infty,0)}(p_{\lambda}(x)
\mbox{ has at least $k$
explosions})\\&=& \PP_{(\infty,-\lambda)}(p_0(x) \mbox{ has at
least $k$ explosions}) \\&=&
\int_{\RR^k}\kappa(-\lambda,dx_1)\kappa(x_1,dx_2)\ldots
\kappa(x_{k-1},dx_k).
\end{eqnarray*}
The second equality uses  the obvious translation
equivariance of $p$.
\end{proof}

In addition:

\begin{proposition}
\label{ToInfty}
 The
$\Hb$ eigenvalues are a.s.\ distinct  with no accumulation point.
Especially,
$\Lambda_k\to \infty$ as $k \to \infty$ a.s., and the minimization
 procedure
of the previous section exhausts the full $\Hb$ eigenvalue set.
\end{proposition}

\begin{proof}
The distinctness follows from the fact that $\PP( N(\lambda+) - N(\lambda) \le 1) = 1$.
Next, starting
at any time $x$, there is positive probability
$\kappa(x,\{\infty\})$ that $p_0$ will converge to $\infty$
without blowing up to $-\infty$. Monotonicity implies that
$\kappa(x,\{\infty\})$ is increasing in $x$. Thus the
number of eigenvalues below $\lambda$ is dominated by a
geometric random variable with parameter
$\kappa(-\lambda,\{\infty\})$, whence it is finite. The
claim follows.
\end{proof}

We mention that second half of Proposition \ref{ToInfty} may also be established by making sense
of the resolvent operator of $\Hb$ and showing it maps the unit ball in $L^*$ onto
H\"older$(3/2)^-$ functions vanishing at infinity.  The simplicity of the above proof
demonstrates the advantages of the diffusion picture.

This connection between the limiting top eigenvalues of the
random matrix ensembles and the explosion time of a simple,
one-dimensional diffusion is new even in the deeply studied
cases of $\b = 1,2$ or $4$.  We now recall the formulas of
Tracy-Widom which,  in conjunction with our result, produce
the identities
\begin{equation}
\label{TW}
  \PP_{(\infty,\lambda)} \big(  \m_1= \infty \big) =
    \left\{ \begin{array}{ll}
       \exp \Bigl( -  \frac{1}{2} \int_{\lambda}^{\infty} (s - \lambda) u^2(s) \, ds  \Bigr)
                                          \exp \Bigl( -  \frac{1}{2} \int_{\lambda}^{\infty}  u(s) \, ds  \Bigr), &  \beta = 1, \\
    \exp \Bigl(- \int_{\lambda}^{\infty} (s - \lambda) u^2(s) \, ds \Bigr), &  \beta = 2, \\
                                           \exp \Bigl(- \frac{1}{2}  \int_{ \lambda^{\prime}}^{\infty} (s - \lambda^{\prime}) u^2(s) \, ds  \Bigr)
                                                    \cosh(   \int_{\lambda^{\prime}}^{\infty}  u(s) \, ds  ),
                                                    & \beta=4.    \end{array} \right.
\end{equation}
Here,  $u(s)$ is the solution of $u^{\prime \prime} = s u +
2 u^3$ (Painleve II) subject to $u(s) \sim \Ai(s) $ as $ s
\rightarrow + \infty$, and $\lambda^{\prime} = 2^{2/3}
\lambda$ in the $\b = 4$ distribution. An important problem
for the future is to obtain the equivalent of the closed
Tracy-Widom formulas for general $\beta$.  Even a direct
verification of (\ref{TW}) would be interesting.

\begin{remark} The diffusion \eqref{e.Riccatidiffusion} seems
efficient for simulating Tracy-Widom distributions as well
as distributions of higher eigenvalues. First note that for
$x \ll 1$, $p(x)$ comes down from $+\infty$ like $1/x$.
Also, the more time accumulated inside the parabola
$\rho_{\lambda} = \{ p^2(x) \le  x + \lambda \}$, the less
likely explosion becomes. That is, the typical path which
hits $-\infty$ does so by tunneling through the narrow part
of $\rho_{\lambda}$. In line with these heuristics good simulations of
the general TW$_\b$ distributions may be obtained by
tracking the explosion probability for $p(x)$ begun at say
$p(0) = O(10^3)$ and run only for $O(1)$ time.
\end{remark}

Moving to applications of the Riccati map to the eigenfunctions,
a detailed but standard analysis of the diffusion
\eqref{e.Riccatidiffusion} using domination arguments shows
that for each $\lambda$, with probability one,
$p(x)/\sqrt{x}\to 1$ (after a finite number of initial
blowups and restarts).  By Fubini, this holds for almost
all $\lambda$. By monotonicity, this happens for all
$\lambda$ except for eigenvalues.  Thus we get

\begin{proposition}\label{p.supexp}
A.s.\  for all solutions $f\in \Hloc$ of $\Hb f=\lambda f$,
with $f(0)=0$ we have the following. If  $f$ is not an
eigenfunctions then $f'(x)/(f(x)\sqrt{x})\to 1$. In
particular, if $f$ grows slower than $\exp((2/3-\eps)
x^{3/2})$ then $f$ is an eigenfunction.
\end{proposition}

\begin{remark}\label{r.boundary} Note that the results of this section, and
indeed, the entire paper, easily extend to general initial
boundary conditions.
\end{remark}

We conclude this section with a decay bound on the $\Hb$ eigenfunctions.  Compare
this to the noiseless  ($\beta = \infty$) limit in which case all eigenfunctions are simply
shifts of the Airy function $\Ai(\cdot)$.

\begin{proposition}
Let $f$ be any eigenfunction of $\Hb$. Then, for any $\eps>0$
there is a random constant $C_{f,\eps}$ so that
$$|f(x)| \le C_{f,\eps}\exp(-(2/3-\eps) x^{3/2}).$$
with probability one.
\end{proposition}

\begin{proof}
 Let $p, q$ be the solutions of $\eqref{Riccati}$
corresponding to Dirichlet and Neumann boundary conditions
at $0$ (i.e.\ $q(0)=0$),  with $p = f'/f$ for the
specified eigenfunction $f$. Then $q$ cannot correspond to
an eigenfunction of the Neumann problem and so by
Proposition \ref{p.supexp} (in conjunction with Remark
\ref{r.boundary}) $q(x)/\sqrt{x}\to 1$. From the
differential equation \eqref{Riccati} we have
$$
\frac{d}{dx}(q-p)=-(q-p)(q+p),
$$
and so
\begin{equation}
\label{eig_dif}
(q-p)(x)=C\exp\left(\int_\m^x -
(q+p)(y)\,dy\right),
\end{equation}
for $\m < \infty$ some random time past the final explosion
of $q$ with $p(\m)$ finite.  With the notation
$Q(x)=\int_\m^x q(y)\,dy$, $P(x)=\int_\m^x p(y)\,dy$, and
$R(x)=Q(x)-P(x)$ the above reads
\begin{equation}\label{e.Rode}
R'(x)=C\exp(R(x)-2Q(x)).
\end{equation}
Now $Q(x)= 2/3x^{3/2}(1+o(1))$ and $P(x)=C+\log
|f(x)|\to-\infty$, implying
\begin{equation}\label{e.Rgrow}
R(x)-x^{3/2}\to \infty.
\end{equation}

To finish, it suffices to show that for all $\eps>0$ and
$x>x_0(\eps)$ we have $R(x)\ge 4/3(1-\eps)x^{3/2}$. Assume
the contrary. Then we can find $x_0$ large so that the
right hand side of \eqref{e.Rode} at $x=x_0$ is at most 1
and for $x>x_0$ we have $q(x)>1/2$.
We claim that the solutions of the ODE \eqref{e.Rode}
started  $x=x_0$ are dominated by the solutions of $\hat
R'(x)=1$.  Indeed, when $R(x)\le \hat R(x)$
\begin{eqnarray*}
C\exp(R(x)-2Q(x))&\le& C\exp(\hat R(x)-2Q(x)) \\ & =&
C\exp\left(R(x_0)-2Q(x_0)+\int_{x_0}^x 1- 2q(y)\,dy
\right)\le 1
\end{eqnarray*}
so that the monotonicity can be maintained. Thus $R(x)\le
\hat R(x)=C+x$ for all large $x$, contradicting
\eqref{e.Rgrow}.
\end{proof}

\section{Tracy-Widom tail bounds: an application of SAE$_\beta$}

This section contains the proof of Theorem \ref{shapethm}.

With $TW_{\beta} = - \Lambda_0(\beta)$ now defined by
(\ref{variationalProbH}), both the upper bound on
$\PP(TW_{\beta} < -a )$ and the lower bound on $P(TW_{\beta}
> a)$  follow from suitable choices of test function $f$ in
$$
    \fip{ f,{\mathcal H}_{\b} f}  = \int_0^{\infty} [  f^{\prime}(x)^2
    + x f^2(x) ] dx
  +  \frac{2}{\sqrt{\beta}} \int_0^{\infty} f^2(x) db_x \ge \Lambda_0
     \, \| f\|_2^2.
$$
The other two bounds run through the Riccati
correspondence.

\bigskip

\noindent {\bf Lower bound, right tail. } Begin with the
observation that
\begin{eqnarray}
\label{test1}
  \PP( TW_{\b} > a ) =
  \PP( \Lambda_0(\b) < -a) & \ge &  \PP \Bigl( \fip{ f,{\mathcal H}_{\b} f } <  - a   \nip{ f,f}  \Bigr)    \\
               &  = &  \PP \Bigl(   \frac{2}{\sqrt{\b}}  \| f\|_4^2 \,
  \mathfrak{ g} < - a \| f\|_2^2 - \| f^{\prime}\|_2^2 - \| \sqrt{x} f\|_2^2 \Bigr)
\nonumber
\end{eqnarray}
for any choice of $f \in L^*$, and  $\mathfrak{g}$  a
standard Gaussian variable. Here we have used the rule that
for $h$ deterministic,  $\int_0^{\infty} h db$ is a
centered Gaussian with variance $\|h\|_2^2$.

Wishing to maximize this probability, the observation is:
for there to be a large negative eigenvalue the random
potential, and then so also the eigenfunction $f$, should
be localized.  This leads one to neglect the
$\int_0^{\infty} x f^2$ term and look for the maximizer $f$
of the expression entering the Gaussian tail, namely
 $$\frac{\int_0^{\infty} (a f^2 + f^{\prime\,2})}{(\int_0^{\infty} f^4 )^{1/2}}.$$

Viewing this problem on the whole line and neglecting
boundary conditions, the optimizers can be computed
exactly: $f(x) = c_1\sech(\sqrt{a}\,x+c_2)$.

Note that $\int_{-\infty}^{\infty} \sech^2(x) dx = 2,
  \int_{-\infty}^{\infty} (\sech^{\prime}(x))^2 dx$ $= 2/3,$ and
$ \int_{-\infty}^{\infty}  \sech^4(x) dx = 4/3$.  Let $f(x)
= \sech(\sqrt{a} (x - 1) )$. Then, on $\RR^+$, with $\sim$
denoting asymptotics as $ a \uparrow \infty$ we have
$$
  a\| f\|_2^2  \sim 2\sqrt{a},\qquad
 \|f'\|_2^2 \sim \frac{2}{3} \sqrt{a}, \qquad
 \| \sqrt{x} f \|_2^2 = O (\frac{1}{\sqrt{a}}), \qquad
 \| f \|_4^4 \sim \frac{4}{3\sqrt{a}}.
$$
Further, while $f(0) \neq 0$, it decays exponentially there
as $a\to\infty$, allowing for an admissible modification
which shares the above evaluations. Returning to
(\ref{test1}) we find that
\begin{eqnarray*}
   \PP ( TW_{\b} > a )
       & \ge & \PP  \left(     \frac{2}{\sqrt{\b}} \times  \frac{2 }{\sqrt{3} }  a^{-1/4}  \,  \mathfrak{g}    <  -   {a}^{1/2} \left(2+\frac{2}{3} +o(1)\right)  \right),
\end{eqnarray*}
producing the desired bound from the simple Gaussian tail
estimate $\PP ( \mathfrak{g} > c) = e^{-c^2(1/2+o(1))}$.
\bigskip

\noindent {\bf Upper bound, left tail. } The reasoning is
the same as that just employed, though in minimizing the
right hand side of $\PP( TW_{\b} < -a) \le \PP ( \fip{f,
{\mathcal H}_{\b} f } >  a \nip{f,f})$ one expects the
optimal $f$ to be relatively ``flat". Neglecting the $\int
(f^{\prime})^2$ term leads to the choice
$$
f(x) =   (x\sqrt{a} ) \wedge  \sqrt{(a - x)^+}\wedge
(a-x)^+.
$$
The middle term is dominant, while the others control
 $\|f'\|_2$. Then
$$
  a\|f\|_2^2 \sim  \frac{a^3}{2}, \qquad
  \| f'\|_2^2 = O(a), \qquad
  \| \sqrt{x} f \|_2^2  \sim  \frac{a^3}{6} , \qquad
   \| f \|_4^4 \sim  \frac{a^3}{3}  .
$$
The proof is completed by substitution
\begin{eqnarray*}
  \PP ( TW_{\b} < -a )
   &  \le   &  \PP  \left( \frac{2}{\sqrt{\b}}   \times \frac{1}{\sqrt{3}} \, {a}^{3/2} \, \mathfrak{g}>
   a^3\left(\frac{1}{2}-\frac{1}{6}+o(1)\right) \right)
  = \exp \Bigl( - \frac{\b}{24} \,a^3 (1+o(1)) \Bigr).
\end{eqnarray*}

\bigskip

\noindent {\bf Lower bound, left tail. } For this we use
the diffusion description of Theorem \ref{difthm}, namely
$\PP( TW_{\b} < - a) = \PP_{(\infty, -a)} ( p \mbox{ never
explodes} )$, where $p$ is the diffusion \eqref{diffusion},
and the subscript indicates the space-time starting point.
By monotonicity,
\begin{eqnarray*}
\lefteqn{
  \PP_{(\infty,-a)} \Bigl( p \mbox{ never explodes} \Bigr)  \ge   \PP_{(1,-a)}
  \Bigl(  p \mbox{ never explodes}  \Bigr) } \\
  &  &  \ge  \PP_{(0,-a)} \Bigl( p(x) \in [0, 2] \mbox{ for all } x\in[-a,0] \Bigr)
 \PP_{(0,0)} \Bigl(  p \mbox{ never explodes} \Bigr).
\end{eqnarray*}
The last factor in line two is some positive number not
depending on $a$. To bound the first probability from
below, employ the Cameron-Martin-Girsanov formula
\begin{eqnarray*}
\lefteqn{
  \PP_{(1,-a)} \Bigl(  p(x) \in [0, 2] \mbox{ for all } x \in[-a,0]  \Bigr) } \\
   &  =  &
   \EE_{(1,-a)} \Bigl[\exp\left( - \frac{\b}{4} \int_{-a}^0 (x - b^2_x) db_x - \frac{ \b}{8} \int_{-a}^0 (x - b^2_x)^2 dx \right);  \, b_x \in [0,2] \mbox{ for all } x\le
   0\Bigr],
\end{eqnarray*}
where $b_x$ is a Brownian motion with diffusion coefficient
$2/\sqrt{\b}$.  On the event in question,
$$
   \frac{\b}{8} \int_{-a}^{0} (x - b^2_x)^2 \,dx  = \frac{ \b}{24} a^3  + O(a^2)   ,
$$
and
$$
    \int_{-a}^0 (x  - b^2_x) db_x =  a b_{-a} + \frac{1}{3} (b^3_{-a} - b^3_0) + (\frac{4}{\beta} - 1)  \int_{-a}^0 b_x dx= O(a).
$$
To finish, note that $\PP_{(-a,0)} ( b_x \in [0,2] \mbox{ for
} x \le 0 ) \ge e^{-ca}$, and so does not interfere with
the leading order.

\bigskip
\noindent {\bf Upper bound, right tail. }  We thank Laure
Dumaz for helping correct an earlier, flawed attempt. The
present proof is based on her master's thesis \cite{Du},
which contains more precise estimates. Write, for $a \gg 1$
and a large $c>0$ to be determined,
\begin{equation}
\label{chain}
  \PP (TW_{\beta} > a ) = \PP_{\infty}(  \mathfrak{m}_{-\infty}< \infty ) \le \PP_{\sqrt{a}-c} (\mathfrak{m}_{-\sqrt{a}} < \infty) ,
\end{equation}
where $\mathfrak{m}$ denotes  the passage time to the indicated level of
the process
$$
   dp(x) = \frac{2}{\sqrt{\beta}} db_x + (a + x - p^2(x)) dx.
$$
To bound the rightmost probability in (\ref{chain}) we introduce a further control on the paths,
and show that there is a numerical constant $c'$ so that
\begin{equation}
\label{pathcontrol}
     \PP_{\sqrt{a}-1} (\mathfrak{m}_{-\sqrt{a}} < \infty) \le  c' \PP_{\sqrt{a}-1} \Bigl( \mathfrak{m}_{-\sqrt{a}} < \infty, \mathfrak{m}_{\sqrt{a}} >  \mathfrak{m}_{-\sqrt{a}} \Bigr).
\end{equation}
This is accomplished by two applications of the (strong) Markov property. From now on we denote
$\mathfrak{m}_{\pm} =  \mathfrak{m}_{\pm \sqrt{a}}$ and
$\mathcal{A} = \{  \mathfrak{m}_{+} >  \mathfrak{m}_{-} \}$, and have
\begin{eqnarray*}
       \PP_{\sqrt{a}-c} (\mathfrak{m}_{-} < \infty) & = & \PP_{\sqrt{a}-c} \Bigl( \mathfrak{m}_{-} < \infty, \mathcal{A} \Bigr) +
          \PP_{\sqrt{a}-c} \Bigl( \mathfrak{m}_{-} < \infty, \mathcal{A}^c \Bigr)  \\
             & \le & \PP_{\sqrt{a}-c} \Bigl( \mathfrak{m}_{-} < \infty, \mathcal{A} \Bigr) +
          \PP_{\sqrt{a}} \Bigl( \mathfrak{m}_{-} < \infty  \Bigr) \\
          & \le &  \PP_{\sqrt{a}-c} \Bigl( \mathfrak{m}_{-} < \infty, \mathcal{A} \Bigr) +
            \PP_{\sqrt{a}} (  \mathfrak{m}_{\sqrt{a}-c}  < \infty)   \PP_{\sqrt{a}-c} \Bigl( \mathfrak{m}_{-} < \infty  \Bigr).
\end{eqnarray*}
Both inequalities use the fact that the hitting probability
of any level below the starting place is decreasing in the
staring time of the diffusion $p$. The desired bound
(\ref{pathcontrol}) then lies in the following claim, the
proof of which we defer to the end.

\begin{claim}
\label{claimconst} There exists a large enough $c$ so that
$\PP_{\sqrt{a}} ( \mathfrak{m}_{\sqrt{a}-c} = \infty)$ is
uniformly bounded below (i.e., independently of $a\gg c$).
\end{claim}

We proceed by performing a change of measure,
\begin{equation}
\label{measurechange}
     \PP_{\sqrt{a}-c} \Bigl( \mathfrak{m}_{-} < \infty, \mathcal{A} \Bigr) =  \lim_{L \ra \infty}  \EE_{\sqrt{a}-c} \Bigl[  R(q, L) \, , \mathfrak{m}_{+} < \mathfrak{m}_{-} < L   \Bigr],
\end{equation}
where $q$ is the diffusion with reversed drift,
$$
   dq(x) = \frac{2}{\sqrt{\beta}} db(x)  + (q^2(x) - x -a) dx,
$$
and the Cameron-Martin-Girsanov factor is given by:
$$
   \log R(q, L) = \frac{\beta}{2} \int_0^{L \wedge \mathfrak{m}_{-}} (a + x - q^2(x) ) dq(x).
$$
An application of It\^o's lemma shows that, for any $z > 0$,
\begin{eqnarray}
\label{expandito}
    \int_0^z ( a + x - q^2(x) ) dq(x)   & =  & a ( q(z) - q(0) ) -  \frac{1}{3} ( q^3(z) - q^3(0) )  \\
                                                             &     &  + \, z q(z) + (4/\beta -1) \int_0^z q(x) dx. \nonumber
\end{eqnarray}
For $z\le \mathfrak{m}_{-}\wedge\mathfrak{m}_{+}$, then
$|q(x)| \le \sqrt{a}$ for $x \in [0,z]$ so the first line
is bounded by a function of $a$, and the second line is
bounded by $(1+|4/\beta -1|) \sqrt{a} \mathfrak{m}_{-} = c_1 \sqrt{a} \mathfrak{m}_{-}$. We
will show that
\begin{equation}
\label{dominated}\log \EE_{\sqrt{a}-c} [ e^{c_1
\sqrt{a}\bar{\mathfrak{m}}},\mathfrak{m}_{+} <
\mathfrak{m}_{-} ] = o(a^{3/2}).
\end{equation}
Letting $L\to \infty$ in \eqref{measurechange}, the
dominated convergence theorem combined with
\eqref{expandito} gives
\begin{equation}\label{laure}
  \PP (TW_{\beta} > a) \le c' e^{-\frac{2}{3} \beta
  a^{3/2}(1+o(1))},
\end{equation}
since for $z = \mathfrak{m}_{-}$ and $q(0) = \sqrt{a} -c$
the first line of (\ref{expandito}) equals $-(4/3) a^{3/2}
+ O(a)$.

Next recall that the expectation \eqref{dominated} is in
terms of the $q$-diffusion, and notice that  on the event
$\{ \mathfrak{m}_{+} < \mathfrak{m}_{-} \}$ the $q$-drift
is bounded above by $-x$.  Therefore, if we introduce the
process  $\bar{q}(x)  = \stb b(x) - x^2/2$ and let $
   \bar{\mathfrak{m}} = \inf \{ x> 0 : \bar{q}(x) = -\sqrt{a} \},
$
it is enough to show that $\log \EE_{\sqrt{a}} [ e^{c \sqrt{a} \bar{\mathfrak{m}}} ] = o(a^{3/2})$ for any constant $c$.  For this, the simple bound
$$
 \PP_{\sqrt{a}} ( \bar{\mathfrak{m}} > t ) \le \PP_{\sqrt{a}} \left(  \stb b(t) > t^2 - \sqrt{a} \right) \le e^{- \frac{\beta}{8} t^3}   \mbox{ for } t > 4 (1 +\stb) \sqrt{a}
$$
will do the job.

To finish, we return to the proof behind the key bound (\ref{pathcontrol}).

\begin{proof}[Proof of Claim \ref{claimconst}] Certainly, the probability that $p$, begun at $\sqrt{a}$, never reaches $\sqrt{a} - c$
is bounded below by the same probability for $p$ replaced by its reflected (downward) at $\sqrt{a}$ version.  Further, when restricted to the interval $[\sqrt{a}-c, \sqrt{a}]$, the $p$-diffusion has
drift everywhere bounded  below by $x$.  Thus we may consider instead the same probability for the appropriate reflected Brownian motion with quadratic drift.

To formalize this, it is convenient to shift orientation.   Let now $\bar{p}(x) = \stb b(x) -x^2/2$, and let $p^*(x)$ denote $\bar{p}$ reflected
(upward) at the origin. Namely,
\begin{equation}
\label{reflected}
    p^{*}(x) = \bar{p}(x) - \inf_{y<x} \bar{p}(y).
\end{equation}
If we can show that  for a large enough $c$,
$
  \PP_0 (  p^*(x) \mbox{ never reaches } c )  > 0,
$
then the $p$-probability in question will also be bounded below, independent of $a$.

What we actually prove is that $M = \sup_{x > 0} p^*(x) $ is a.s.\ finite.  Let $\kappa > 1$, and for each $n \ge 0$ define the event
$$
   D_n = \Bigl\{ \bar{p}(x) \mbox{ hits } (1-n)\kappa \mbox{ for some } x  \mbox{ between } \mathfrak{m}_{-n\kappa} \mbox{ and } \mathfrak{m}_{-(n+1)\kappa}  \Bigr\}.
$$
>From the representation (\ref{reflected}) one sees that $M \ge 2\kappa$ implies that some $D_n$ must occur. Define as well
the event
$$
  A_{\kappa} =  \Bigl\{ \bar{p}(x) \ge -\frac{1}{2} (\kappa x^2 + 1) \mbox{ for all } x > 0 \Bigr\},
$$
and note that, on $A_\kappa$, $\mathfrak{m}_{-n\kappa} \ge \sqrt{ 2n -1}$.   Hence, still on $A_\kappa$, for all $n \ge 1$ the shifted process
 $\bar{p}(\mathfrak{m}_{-n\kappa} + x)  - \bar{p}(\mathfrak{m}_{-n\kappa})$ is dominated by the scaled Brownian motion $x \mapsto
 \stb b(x)$ plus constant drift $ -(1/2) \sqrt{2n-1}$.
It follows that: for $n \ge 1$,
\begin{eqnarray}
\label{constdrift}
   \PP ( D_n  \cap A_\kappa | \mathcal{F}_{\mathfrak{m}_{-n\kappa}} ) & \le & \PP_0 \Bigl( x \mapsto  ( \stb b(x) -\mbox{$\frac{1}{2}$}\sqrt{2n-1}\;x ) \mbox{ hits } \kappa \mbox{ before } -\kappa
     \Bigr)   \\
   & =  &   \frac{1}{ 1 + e^{ \kappa (\beta/2) \sqrt{2n-1} }}. \nonumber
   \end{eqnarray}

Putting the above together we have
\begin{eqnarray*}
   \PP (M > 2\kappa) & \le & \PP (A_\kappa^c ) + \PP ( A_\kappa \cap ( \cup_{n \ge 0} D_n )  )
   \le \PP(A_\kappa^c ) + \PP (D_0) +\sum_{n\ge 1}\PP( A_\kappa \cap D_n
   ).
\end{eqnarray*}
The sum of the series converges to $0$ by
\eqref{constdrift}. Since $D_0$ implies that $\sup_{x>0}
\bar{p}(x)
> \kappa$, both $\PP (A_\kappa^c) $ and $\PP(D_0) $ tend to
zero as $\kappa \ra \infty$ by the Law of the Iterated
Logarithm.
\end{proof}

\section{Convergence of discrete models and universality}

This section establishes a general and rather weak set of
conditions under which the bottom eigenvalues of random
symmetric tridiagonal matrices converge to the bottom
eigenvalues of a corresponding stochastic differential
operator.   In many ways this is the central result of the
paper.

To explain the setup, consider a sequence of discrete-time
$\mathbb R^2$-valued random sequences
$((y_{n,1,k},y_{n,2,k}); 1\le k \le n)$. Let $m_n=o(n)$ be
a scaling parameter.  (In the particular case of the
Hermite  $\beta$-ensembles we have $m_n=n^{1/3}$.) For each
$n$, we build an $n\times n$ tridiagonal matrix $H_n$.

Let $T_n$ denote the shift operator $(T_n v)_k=v_{k+1}$
acting on $\RR_1\times \RR_2\ldots$. Let $(T_n^t
v)_k=v_{k-1}1_{k\ge 1}$ be its adjoint, and let $R_n$ denote
the restriction operator $(R_n v)_k=v_k\one_{k\le n}$.  Let also
$\Delta_n=m_n(I-T_n^t)$ be the difference quotient
operator, and finally set
\begin{equation} \label{H_n}
H_n = R_n\left(-\Delta_n\Delta_n^t + (\Delta_n
y_{n,1})_\times + (\Delta_n y_{n,2})_\times\frac12(T_n +
T_n^t)\right),
\end{equation}
 where the subscript $\times$
denotes element-wise multiplication by the corresponding
vector. Then $H_n$ maps the coordinate subspace $\RR^n\to
\RR^n$, and its matrix with respect to the coordinate basis
in $\RR^n$ is symmetric tridiagonal with $(2m_n^2+m_n
(y_{n,1,k}-y_{n,1,k-1}),k\ge 1)$ on the diagonal and
$(-m_n^2+ m_n(y_{n,2,k}-y_{n,2,k-1})/2, k \ge 1)$ below and
above the diagonal.  Roughly speaking, $H_n$ is the
discrete Laplacian plus integrated potential
$y_{n,1}+y_{n,2}$.

Additionally, define $y_{n,i}(x)=y_{n,i,\lfloor x m_n\rfloor
}1_{xm_n\in[0,n]}$. By choice, $(\sqrt{m_n}\times
y_{n,i,k},k\ge 0)$ is on the scale of a simple random walk,
so no additional spatial scaling will be required.

Our basic convergence result rests on two sets of
assumptions on the processes $y_{n,i=1,2}$.

\bigskip

\noindent {\em Assumption 1 (Tightness/Convergence) } There
exists a continuous process $x \mapsto y(x)$ such that
\begin{eqnarray}\nonumber
\big(y_{n,i}(x);\;x\ge 0\big) &&i=1,2\quad \mbox { are tight in law, }\\
\big(y_{n,1}(x)+y_{n,2}(x);\; x\ge 0\big) &\Rightarrow
&\big(y(x);\, x\ge 0\big) \quad \mbox{ in
law,}\label{condition1}
\end{eqnarray}
with respect to the Skorokhod topology of paths, see
\cite{EthierKurtz} for definitions.

\bigskip

\noindent {\em Assumption 2 (Growth/Oscillation bound) }
There is a decomposition
 \begin{equation}\label{e.assumptiongrowth}
y_{n,i,k} =m_n^{-1}\sum_{\ell=1}^{k} \eta_{n,i,\ell} \, +
\, w_{n,i,k}
 \end{equation} for $\eta_{n,i,k}\ge 0$,
 deterministic,
unbounded nondecreasing continuous functions
$\bareta(x)>0,\zeta(x) \ge 1$, and random constants
$\kappa_n(\omega)\ge 1$ defined on the same probability
space which satisfy the following. The $\kappa_n$ are tight
in distribution, and, almost surely,
\begin{eqnarray}
 \bareta(x)/\kappa_n\;-\kappa_n\le\; \eta_{n,1}(x)  + \eta_{n,2}(x) &\le& \;\kappa_n(1+\bareta(x)),
\label{bounds1}
 \\
\label{bounds2} \eta_{n,2}(x)&\le &2m_n^2\\
 |w_{n,1}(\xi)-w_{n,1}(x)|^2  +  |w_{n,2}(\xi)-w_{n,2}(x)|^2   & \le& \kappa_n(1+\bareta(x)/\zeta(x)).
\label{bounds3}
\end{eqnarray}
for all $n$ and $x,\xi\in[0,n/m_n]$ with $|x-\xi|\le 1$.

\bigskip

We may now define the limiting operator. Just as in Section
\ref{s.basic} we note that
\begin{equation}\label{e.generalH}
 H=-\frac{d^2}{dx^2} +y'(x)
\end{equation}
 maps $\Hloc$ to the space $D$ of distributions
via integration by parts. Without changing the notation, we
generalize the Hilbert space $L^*\subset L^2(\RR^+)$
introduced there. This consists of functions with $f(0)=0$
and
$$
\|f\|_*^2=\int_0^\infty f'(x)^2 + (1+\bareta(x))f^2(x)
dx<\infty.
$$
The eigenvalues and eigenfunctions are defined again as
$(\lambda,f)\in\RR\times  L^*\setminus\{0\}$ with $\| f \|_2 =1$ satisfying
\eqref{e.generalH}. Recall from Section \ref{s.basic} that
this means
$$
f'(x)=\int_0^x -y(z)f'(z) -\lambda f\, dz + f(x)y(x),
$$
or, equivalently, for all $\varphi \in C_0^\infty$ it holds
$$
\int f \varphi'' dx = \int -\lambda f \varphi -y
f'\varphi-yf\varphi' \,dx.
$$

With the picture laid out, a few words on Assumptions 1 and
2 are in order. The former simply asks for correspondence
between $H_n$ and $H$ at the level of integrated
potentials.  The latter, more technical condition, will
imply the compactness necessary to maintain discrete
spectrum as $n \uparrow \infty$.

\begin{theorem}[Convergence in law] \label{weak}
Given Assumption 1 and 2 above and any fixed $k$,
 the bottom $k$ eigenvalues of the matrices $H_n$ converge in
law to the bottom $k$ eigenvalues of the operator $H$.
\end{theorem}

We will also show that, after a natural embedding, the
eigenfunctions also converge in $L^2$.

\begin{proof}[Proof: Reduction to the deterministic case] It will
be convenient to find subsequences along which we have
limits for all desired quantities.

The upper bound (\ref{bounds1}) shows that $(\int_0^x
\eta_{n,i}(t)dt;\;x\ge 0)$ is tight in distribution for
$i=1,2$. For any subsequence we can extract a further
subsequence so that we have joint distributional
convergence
\begin{eqnarray} \nonumber
 (\int_0^x \eta_{n,i}(t)dt;\;x\ge 0)
 &\Rightarrow& (
 \eta^\dagger_{i}(t)dt;\;x\ge 0), \\
( y_{n,i}(x);\;x\ge 0) &\Rightarrow& (y_i(x);\;x\ge 0),
\label{weakconv}
\\
\kappa_n&\Rightarrow & \kappa, \nonumber
\end{eqnarray}
where the convergence in the first lines is in the
uniform-on-compacts topology, and the second, in the
Skorokhod topology. Then by Skorokhod's representation
theorem (see Theorem 1.8, Chapter 2 of \cite{EthierKurtz}) we can realize
this convergence as a.s.\ convergence on some probability
space so that the conditions of Proposition \ref{deter}
below are satisfied with probability one.

Note that (\ref{bounds1}-\ref{bounds2}) are local Lipschitz
bounds on the $\int \eta_{n,i}$, and so they are inherited
by their limit  $\eta^\dagger_{i}$. Thus
$\eta_i=(\eta^\dagger_i)'$ is defined  almost everywhere,
and satisfies (\ref{bounds1}-\ref{bounds2}). Further,
$\eta_i$ can be defined everywhere so that
(\ref{bounds1}-\ref{bounds2}) continues to hold.

It also follows that each $w_{n,i}=y_{n,i}-\sum \eta_{n,i}$
must have a limit which we denote $w_i$. We further denote
$w=w_1+w_2$ and $\eta=\eta_1+\eta_2$. The claim now follows
from Proposition \ref{deter} below.
\end{proof}

\begin{proposition}[Deterministic convergence]\label{deter}
Assume that the each of the convergence statements in
(\ref{weakconv}) hold deterministically and that  the
bounds (\ref{bounds1}-\ref{bounds3}) hold with some
deterministic constant $\kappa$ replacing $\kappa_n$. Then,
for any $k$, the lowest $k$ eigenvalues of the matrices
$H_n$ converge to the lowest $k$ eigenvalues of $H$.
\end{proposition}

In the next subsection, we establish some properties of the
limiting operator. Afterwards, we prove Proposition
\ref{deter}.

\subsection*{Properties of the limiting operator}

Just as in Section \ref{s.basic}, we extend the bilinear
form $\fip{\cdot,H \cdot}$ from $C_0^\infty \times L^*$ to
$L^*\times L^*$. We want to define the extension as
$\fip{f,H f}:=\int f'^2 +\eta f^2 dx+ \int f^2 d w$ (with
$\fip{f,Hg}$ then defined by polarization), but we still
need to define and control the last term. By the next
lemma, this can be done via the integration by parts
already employed in \eqref{e.quadraticform}:
\begin{equation}
\int_0^{\infty} f^2(x)\, dw_x = \int_0^{\infty} f^2(x)\,
(w_{x+1} - w_x)\,dx + 2\int_0^{\infty} f'(x) f(x)
(\int_x^{x+1} w_t dt - w_x)\, dx.
\label{eq:DecompStochTerm1}
\end{equation}

\begin{lemma} \label{*extension}
The integrals on the right of (\ref{eq:DecompStochTerm1})
are defined and finite for $f\in L^*$. Moreover there exist
$c_8,c_9,c_{10}>0$ so that
$$
c_8 \|f\|^2_{*}-c_9\|f\|^2_2\le \fip{f,H f}\le c_{10}
\|f\|^2_{*}.
$$
\end{lemma}

\begin{proof}
By taking limits of the inequalities
(\ref{bounds1}-\ref{bounds3}) on $\eta_{i,n}$ and $w_n$ we
get bounds for $\eta_i$ and $w$. In particular
$\max(|w_{x+1}-w_x|,|w_{x+1}-w_x|^2)\le c_\eps+\eps\bareta$
where $\eps$ can be made small. Now we write
$\fip{f,Hf}=A+B$ where $B$ is the fluctuation term
(\ref{eq:DecompStochTerm1}), and the potential term
satisfies  $\frac 1 {\kappa}\|f\|_*^2-c\|f\|_2^2\le A \le
c\|f\|^2_*$. To bound $B$, first write
$$\int_0^{\infty} f^2(x)\, |w_{x+1} - w_x|\,dx \le
 \s{f,(c_\eps+\eps)\bareta) f} .$$
 For the second term, we average the inequality
 $\sup_{|x-\xi|\le 1}|w_{\xi}-w_x|\le
 |c_\eps+\eps \bareta(x)|^{1/2}$ and use an inequality of means:
$$
2\int_0^{\infty} |f'(x) f(x) (\int_x^{x+1} w_t dt - w_x)|\,
dx \le \sqrt{\eps}\|f'\|_2^2 +
\s{f,\frac{1}{\sqrt{\eps}}(c_\eps+\eps\bareta) f}.
$$
The above inequalities give $|B|\le
2\sqrt{\eps}\|f\|^2_*+c'_\eps\|f\|_2^2$. Setting $\eps$
small we get the results.
\end{proof}

The bounds immediately imply the following.

\begin{corollary}\label{continuity} (i) The bilinear form $\fip{\cdot,H \cdot}:L^*\times L^*\to \RR$ is
continuous. (ii) It does not depend on the decomposition
$y=w+\int \eta$. (iii) The eigenvalues and eigenfunctions
$(\lambda,f)$ of $H$ satisfy $\fip{g,Hf}=\lambda\s{g,f}$
for all $g\in L^*$. (iv) In particular,
$\fip{f,Hf}=\lambda\s{f,f}$.
\end{corollary}

\begin{proof}[Proof of the Corollary]
Since $L^*\subset L^2$ is a continuous embedding, it
suffices to prove that $\fip{\cdot,H \cdot}+
\,c_9\,\s{\cdot,\cdot}$ is continuous. This form is
nonnegative definite by the lemma. Continuity {\it (i)} now
follows from Cauchy-Schwarz applied to the form and the
bounds of the Lemma. For
 {\it (ii)} and {\it (iii)} approximate $g$ by smooth compactly supported functions and use
continuity.
\end{proof}

Together, these two statements provide an analogue of Lemma
\ref{Compactness} in a more general context. In particular,
the (discrete) eigenvalues of $H$ may now be defined
variationally. The arguments used in Section \ref{s.basic}
in conjunction with the limiting bounds of
(\ref{bounds1}-\ref{bounds3}) give the following.

\begin{lemma} The lowest $k$ elements of the set of
eigenvalues of $H$ exist  and admit the variational
characterization via the bilinear form $\fip{f,H f}$.
\end{lemma}

\subsection*{Tightness}
 Next, we define the discrete analogue of the norm
$\|\cdot \|_*$. For $v\in \RR^n$ let
$\|v\|_2^2=m_n^{-1}\sum_{k=1}^n v_k^2$ with scaling to
match the continuum norm. Let
$\bareta_{n,k}=\bareta(k/m_n)$ and let
\begin{equation}
\label{discretenorm} \|v\|_{*n}^2 = \|\Delta_n v \|^2_2 +
\|  (\bareta_{n,\cdot})^{1/2}  v\|^2_2+ \| v\|^2_2,
\end{equation}
and note that $\|v\|_2 \le \|v\|_{*n}$.  We continue with a
bound on $H_n$.

\begin{lemma} \label{compact} Assume the bounds (\ref{bounds1}-\ref{bounds3}).
Then there are constants $c_{11},c_{12},c_{13}>0$ so that
for all $n$ and all $v$ we have
$$
c_{11} \|v\|^2_{*n} - c_{12} \|v\|^2_2\le \langle v,H_nv
\rangle \le c_{13} \|v\|^2_{*n}.
$$
\end{lemma}
\begin{proof}
Drop the subscript $n$, and recall the definition of the
difference quotient $\Delta v_k=m(v_{k+1}-v_{k})$. We
recall the following consequence of the discretized bounds
(\ref{bounds1}-\ref{bounds3}):
\begin{eqnarray*}
\bareta_k/\kappa-\kappa &\le& \eta_{1,k}+\eta_{2,k} \;\le
\kappa\bareta_k+\kappa
\\
\eta_{2,k} &\le& 2m^2,
\\
|w_{i,\ell}-w_{i,k}|^2 &\le&  \eps \bareta_k + c_\eps ,
\qquad k \le \ell \le {k+m}, \qquad {i=1,2}.
\end{eqnarray*}
Here $\eps$ can be arbitrarily small at the expense of
$c_\eps$. Let $w_k=w_{1,k}$, $u_k=w_{2,k}$. By definition
of $H_n$,
 \begin{equation}\label{label} m\langle v,H_n v\rangle = \sum_{k=0}^n \left((\Delta v_k)^2 +
\eta_{2,k} v_k v_{k+1}+\eta_{1,k}v_k^2\right) +
\sum_{k=0}^n\left(\Delta w_k v_k^2 + \Delta u_k
v_kv_{k+1}\right).
 \end{equation}
Let $A$, $B$ denote the two sums. Using the inequality
$ab\ge -(b-a)^2/3+a^2/4$, we obtain the following lower
bound for the second summand in $A$:
\begin{eqnarray}\nonumber
\eta_{2,k} v_k v_{k+1}
  \ge
-\frac{\,\eta_{2,k}}{3}(v_{k+1}-v_k)^2 + \frac{\eta_{2,k}
v_k^2}{4} \ge -\frac{2 }3 (\Delta v_k)^2 +
\frac{\bareta_{2,k} v_k^2}{4\kappa}-\frac{\kappa v_k^2}{4},
\label{bb2}
\end{eqnarray}
thus we get $A\ge \frac {m}{4\kappa}\|v\|_*^2-cm\|v\|^2$.
We also clearly have $|A|\le cm\|v\|_*^2$.
 To bound $|B|$, set $\delta w_k= (w_{k+m}-w_k)$.
Summation by parts gives
\begin{eqnarray}\nonumber
  \sum_{k=0}^n \Dws_k  v_k^2
  &=&
  \sum_{k=0}^n \left( \Dws_k-\dws_k \right) v_k^2
  + \sum_{k=0}^n \dws_k  v_k^2
  \\  \label{firstsummand} &=&
  \sum_{k=0}^n \left(
  \sum_{\ell=k+1}^{k+m}(\ws_\ell-w_k) \right) (v_{k+1}^2-v_k^2)
  + \sum_{k=0}^n \dws_k v_k^2.
\end{eqnarray}
By the bound on $w$ the absolute value of the first summand
in (\ref{firstsummand}) is not more than
$$
m\left(\eps\bareta_k +c_\eps\right)^{1/2}
|v_{k+1}^2-v_{k}^2| \le
\frac{1}{\sqrt\eps}\left(\eps\bareta_k+
c_\eps\right)(v_k+v_{k+1})^2+
\sqrt{\eps}m^2(v_{k+1}-v_{k})^2.
$$
Together with the bound on the second sum in
(\ref{firstsummand}) this yields
\begin{equation}\label{bb} m \left|\sum_{k=0}^n \Dws_k
v_k^2 \right|\le (\frac{1}{\sqrt{\eps}}+1)\sum_{k=0}^n
\left(\eps\bar \eta_k + c_+\eps \right) v_k^2+
\sqrt{\eps}m^2\sum_{k=0}^n (v_{k+1}-v_{k})^2.
\end{equation}
The argument for $u$  starts with the same summation by
parts, with $v_kv_{k+1}$ playing the part of $v_k^2$. After
bounding the $u$ terms, we use the  inequalities
$2|v_kv_{k+1}|\le v_k^2+v_{k+1}^2$ and
$$|v_{k}v_{k-1}-v_kv_{k+1}|\le |v_k\|v_k-v_{k-1}| +
|v_k\|v_{k+1}-v_k|$$ together with Cauchy-Schwarz to get an
estimate of the form (\ref{bb}) for $\sum_{k=0}^n \Delta
u_k v_kv_{k+1}.$ Thus we find that $|B|\le
c\sqrt{\eps}m\|v\|_*^2+c_\eps'm\|v\|_2^2$. For $\eps$
sufficiently small the claims follow.
\end{proof}

\subsection*{Operator convergence}

\noindent {\bf Embedding vector spaces.} We embed the
domain $\RR^n$ of $H_n$ in $L^2(\RR^+)$ in an isometric
way, identifying $v\in \RR^{n}$ with the step function
$v(x)=v_{\lceil m_n x\rceil}$ supported on $[0,n/m_n]$. Let
$L^*_n$ denote the space of such step functions, and let
${\mathcal P}_n$ denote the $L^2$-projection to this space.
Let $(T_n f)(x) = f(x+m_n^{-1})$ denote the shift operator,
and let $R_n(f)=f\one_{[0,n/m_n]}$ denote the restriction.
Let $\Delta_n =m_n(I-T_n^t)$. These operators are simply
extensions of the already defined action of $T_n$ and
$\Delta_n$ on $L^*_n$. Thus the formula (\ref{H_n}) extends
the definition of $H_n$ to $L^2$.

It is easy to check the following: (i) ${\mathcal P}_n$ and
$T_n$, and so $\Delta_n$ commute; (ii) for $f\in L^2$ we
have ${\mathcal P}_nf \to f$ in $L^2$; and (iii) when
$f'\in L^2$ and $f(0)=0$ we have $\Delta_n f \to f'$ in
$L^2$.

\begin{lemma} \label{mainconv}
Assume that $f_n \in L^*_n$ and $f_n\to  f$ weakly in $L^2$
and $\Delta_n f_n \to f'$ weakly in $L^2$. Then for all
 $\varphi \in C_0^\infty$ we have
$\s{\varphi,H_n f_{n}} \to \fip{\varphi,H f}$. In
particular
\begin{equation}\label{e.fHf}
\nip{\mathcal P_n \varphi,H_n \mathcal P_n
\varphi}=\nip{\varphi,H_n \mathcal P_n
\varphi}\to\fip{\varphi, H \varphi}.
\end{equation}
\end{lemma}

\begin{proof} Because we are dealing with $\varphi$ of compact support, may drop the restriction part $R_n$ for $H_n$.
The convergence
$$
\nip{\varphi,\Delta_n\Delta_n^t
f}=\nip{\Delta_n\Delta_n^t\varphi, f} \to
\nip{\varphi'',f}=\fip{\varphi,f''}
$$
is clear, so it remains to check the potential term.  First
note that if $I$ is a finite interval, and $g_n \to_{L^2}
g$ and $h_n \to h$ is $L^2(I)$-bounded and converges weakly
in $L^2(I)$ then
\begin{equation}\label{supl2}
%\left|\int_I gh-g_nh_n\,dx\right| \le \sup_{x\in
%I}\left|g_n(x)-g(x)\right|\int_I |h_n|\,dx+\left|\int_I g
%(h_n-h)\,dx\right| \to 0
\s{g_n,h_n \one_I} \to \s{g,h\one_I}.
\end{equation}

Let $I$ be a finite closed interval supporting
$\Delta_n\varphi$, $\varphi'$ and $\varphi$.  The potential
term is
$$
\s{\varphi, \left((\Delta_n y_{n,1})_{\times} + (\Delta_n
y_{n,2})_{\times}\mbox{$\frac12$}(T_n+T_n^t)\right)f}.
$$
Setting $y_n=y_{n,1}+y_{n,2}$, we first approximate the
right hand side by
\begin{eqnarray*}
\s{\varphi, (\Delta_n y_{n})_{\times}f_n} &=&
 \s{\Delta_n^t (\varphi
f_n) , y_n } \\&=& \s{\varphi \Delta_n^t f_n + f_n
\Delta_n^t\varphi + m_n^{-1}\Delta_n^tf_n \Delta_n^t
\varphi , y_n }
\\&=& \s{\Delta_n^t f_n,  \varphi y_n } + \s{f_n  ,y_n\Delta_n^t\varphi
 } +  m_n^{-1} \s{\Delta_n^tf_n , y_n \Delta_n^t \varphi }.
\end{eqnarray*} The
first two terms in the above converge to the desired limits
by (\ref{supl2}), and the last one converges to 0 because
it is bounded without the extra scaling term. The error
term in the above approximation comes as a sum of  $T_n$
and $T_n^t$ terms; we consider twice the $T_n^t$ part:
\begin{eqnarray}\nonumber
|\s{\varphi, (\Delta_n y_{n,2})_{\times}(I-T_n^t)f_n}| &=&
|\s{\varphi\, m_n^{-1}\Delta_n y_{n,2},\Delta_nf_n}|
\\&\le& \|m_n^{-1}\Delta_n y_{n,2}\one_I\|_2\,
\|\Delta_nf_n\|_2\, \sup_{x\in \RR} |\varphi(x)|.
\label{wiggle}
\end{eqnarray}
Now $m_n^{-1}\Delta_n y_{n,2}\one_I$ is the restriction to
$I$ of the difference $y_{n,2}-T_n^ty_{n,2}$, in which both
terms converge to $y_2$ in the Skorokhod topology. In
particular, they converge a.e., and since they are locally
bounded, their difference converges locally in $L^2$ to  0.
This shows that \eqref{wiggle} vanishes in the limit. We
handle the $T_n$ term similarly.
\end{proof}

\begin{lemma} \label{smooth}  Recall the discrete $\| \cdot \|_{*n}$ norm from
(\ref{discretenorm}). Assume that $f_n\in L^*_n$,
$\|f_n\|_{*n} \le c$, and $\|f_n\|_2=1$. Then there exists
$f\in L^*$  and a subsequence $n_k$ so that
$f_{n_k}\to_{L^2} f$ and for all $\varphi\in C_0^\infty $
we have $\langle \varphi,H_{n_k} f_{n_k} \rangle \to
\fip{\varphi,H f}.$
\end{lemma}

\begin{proof}
Since $f_n$ and $\Delta_nf_n$ are bounded in $L^2$, we can
find a subsequence along which  $f_n\to f\in L^2$ and
$\Delta_n f_n\to \tilde f\in L^2$ weakly. Considering
$\s{\Delta_n f_n,1_{[0,t]}}$ we get that $\int \tilde f=
f$, that is $f$ has a differentiable version and $\tilde f
= f'$. The bounded nature of the $\bar \eta$ terms in the
$L^*_n$-norm gives sufficient tightness so that we have
$f\in L^*$ and $f_n\to_{L^2} f$.  The last part then follows
from Lemma \ref{mainconv}.
\end{proof}

We break up the proof of Proposition \ref{deter} into two
Lemmas. Let $(\lambda_{n,k},v_{n,k}),k \ge 0$ be the lowest
eigenvalues and the embedded normalized eigenfunctions of
$H_n$, and let $(\Lambda_k,f_k)$ be the same for $H$.

\begin{lemma}\label{l.liminf}
For $k\ge 0$ we have $ \underline \lambda_k= \liminf
\lambda_{k,n} \ge \Lambda_k. $
\end{lemma}

\begin{proof}
Assume $\underline \lambda_k<\infty$. Since the eigenvalues
of $H_n$ are uniformly bounded below, we can find a
subsequence so that $(\lambda_{n,1},\ldots,
\lambda_{n,k})\to
(\xi_1,\ldots,\xi_k=\underline\lambda_k)$. By Lemma
\ref{compact}, the corresponding eigenfunctions have $L^*_n$
norm uniformly bounded, and Lemma \ref{smooth} now implies
that for a further subsequence, their $L^2$ limit exists.
Moreover, by the same lemma this limit must consist of
orthonormal eigenfunctions of $H$ with eigenvalues at most
$\underline \lambda_k$. The orthonormality of the eigenfunction set
shows that they correspond to $k$ distinct states and the proof is
finished.
\end{proof}

\begin{lemma} For $k\ge 0$ we have $\lambda_{k,n}\to \Lambda_k$ and
$v_{n,k}\to_{L^2} f_k$.
\end{lemma}

\begin{proof}
For an inductive proof, we assume the claim holds up to
$k-1$. First, we find $f_k^\eps\in C_0^\infty$ $\eps$-close
to $f_k$ in $L^*$. Consider the vector
 \begin{equation}\label{gs}
 f_{n,k}={\mathcal P}_n f_k^\eps-\sum_{\ell=1}^{k-1}\langle v_{n,\ell},{\mathcal P}_n f_k^\eps \rangle
 v_{n,\ell}.
 \end{equation}
 We have a uniform bound on the $L^*_n$ norm of
$v_{n,\ell}$ by Lemma \ref{compact}, and $\|{\mathcal
P}_nf_k^\eps-v_{n,k}\|_2\le \|{\mathcal
P}_nf_k^\eps-f_k^\eps\|_2+\|v_{n,k}-f_k^\eps\|_2$, is, for
large $n$ bounded by $2\eps$. Thus the $L^*_n$-norm of the
sum is bounded by $c\eps$. By the uniform bound $\nip
{v,H_nv}\le c \|v\|_{*n}^2$ of Lemma \ref{compact} and the
variational characterization in finite dimensions we also have that
\begin{eqnarray} \label{standardtrick}\limsup \lambda_{n,k}
 \le
\limsup_{n\to\infty} \frac{\langle f_{n,k},H_n f_{n,k}
\rangle}{\langle f_{n,k}, f_{n,k} \rangle}
 =
  \limsup_{n\to\infty} \frac{\langle {\mathcal P}_n f_k^\eps,H_n {\mathcal P}_n
f_k^\eps \rangle}{\langle {\mathcal P}_n f_k^\eps,
{\mathcal P}_n f_k^\eps \rangle}+o_\eps(1),
\end{eqnarray}
where $o_\eps(1)\to0$ as $\eps\to 0$. Then \eqref{e.fHf} of
Lemma \ref{mainconv} provides
$$
\lim_{n\to \infty} \nip{{\mathcal P}_n
f_k^\eps,H_n{\mathcal P}_nf_k^\eps}= \fip{f_k^\eps,H
f_k^\eps},
$$
and therefore the right hand side of \eqref{standardtrick}
equals
$$ \frac{\fip{ f_k^\eps, Hf_k^\eps
}}{\langle f_k^\eps, f_k^\eps \rangle}+o_\eps(1)
  =  \frac{\fip{ f_k, Hf_k
}}{\langle f_k, f_k \rangle}+o_\eps(1).
$$
Now letting $\eps\to 0$ the right hand side converges to
${\fip{ f_k, Hf_k }}/{\langle f_k, f_k \rangle}=\Lambda_k$.
We have shown $\lambda_{n,k}\to  \Lambda_k$.

Lemma \ref{smooth} implies that any subsequence of the
$v_{n,k}$ has a further subsequence converging in $L^2$ to
some $g\in L^*$ satisfying $Hg=\Lambda_k g$. Thus
 $g=f_k$, and so  $v_{n,k}\to_{L^2} f_k$.
\end{proof}

\section{CLT and tightness for tridiagonal $\b$-ensembles}

\label{checkconditions}

At last we verify that the $\b$-Hermite and Laguerre
ensembles satisfy the conditions
(\ref{condition1})-\eqref{bounds3} of Theorem \ref{weak},
and so complete the proofs of  Theorems \ref{mainthm} and
\ref{Lagthm}.

The following theorem is what we need from the far more
general Theorem 7.4.1 on page 354 in Ethier and Kurtz
\cite{EthierKurtz}. Denote $\Da y_{n,k}=y_{n,k}-y_{n,k-1}$.

\begin{corollary}\label{c.either}
Let $a\in \RR$ and $h\in C_1(\RR^+)$, and let $y_n$ be a
sequence of processes with $y_{n,0}=0$ and independent
increments.  Assume that
$$m_n \EE\Da
y_{n,k} =h'(k/m_n)+o(1), \qquad m_n \EE (\Da
y_{n,k})^2=a^2+o(1),\qquad m_n \EE (\Da y_{n,k})^4 = o(1)$$
uniformly for $k/m_n$ on compact sets as $n\to\infty$. Then
$y_n(t) = y_{n, \lfloor tm_n \rfloor }$ converges in law, with respect to
the Skorokhod topology, to the process $h(t)+ab_t$, where
$b$ is standard Brownian motion.
\end{corollary}

\begin{proof}
The time-homogeneity required in the theorem can be
replaced by introducing a space coordinate recording time.
The supremum increment bound of the theorem follows from
Markov's inequality and the fourth moment bound here.
\end{proof}

\subsection*{The $\beta$-Hermite case}

Starting with the scaled Hermite matrix ensembles $H_n =
{\tilde{H}}_n^{\b}$, we identify $m_n = n^{1/3}$. After
rearranging some terms we find
\begin{eqnarray*}
    {y}_{n, 1,k} &=& w_{n,1,k} \,= -\,n^{-1/6}(2/\beta)^{1/2} \sum_{\ell=1}^{k}
    g_\ell,
\\   y_{n, 2,k}  &=&   n^{-1/6}\sum_{\ell=1}^{k}  2 \Bigl(
\sqrt{n} - \frac{1}{\sqrt{\b } }
 \chi_{\b (n-\ell)} \Bigr).
\end{eqnarray*}
Also, by choosing $ \eta_{n,2,k}=2\sqrt{n}-2\beta^{-1/2}\,\EE
\chi_{\beta(n-k)}$ both $w_{n,1,k}$ and $w_{n,2,k}$ are independent-increment
martingales.
Using the notation and results of Corollary
\ref{c.either} and standard moment computations for the
normal and gamma distributions, we get the following.
\begin{lemma}\label{l.Hermite.conv}
As $n\to\infty$ for the Skorokhod topology we have, in law
$$
y_{n,i}(\cdot) \Rightarrow(2/\beta)^{1/2}b_x+x^2(i-1),
\qquad i=1,2.
$$
\end{lemma}
Independence of the $i=1,2$ cases now implies
\eqref{condition1} of Assumption 1.
\begin{lemma} The bounds \eqref{bounds1}, \eqref{bounds2} of Assumption 2
hold with $\bar{\eta}(x)=x$.
\end{lemma}
\begin{proof}
There is the estimate
%\begin{equation}
% \label{moment0}
$$
   \sqrt{r}(1-4/r)\le \EE \chi_{r}   =  \sqrt{2} \, \frac{  \Gamma( (r+1)/2 )
         }{   \Gamma( r/2)   } \le\sqrt{r},
$$
%\end{equation}
and, with again  $ \eta_{n,2,k}=2\sqrt{n}-2\beta^{-1/2} \EE
\chi_{\beta(n-k)}$,
it follows that
$$
kn^{-1/2}-c\le \eta_{n,2,k}\le 2kn^{-1/2}+c,
$$
where $c$ depends on $\beta$ only.
\end{proof}

Lastly, for \eqref{bounds3} of Assumption 2, it suffices
to prove a tight random constant bound on
$$
\sup_{k=1\ldots n/m_n} k^{\eps-1} \sup_{\ell=0\ldots m_n}
\left\vert w_{n,i,km_n+\ell}-w_{n,i,km_n} \right\vert^2.
$$
(Notice the estimate is being done in blocks.)
Squaring, replacing the first supremum by a sum,  and then taking
expectations gives
$$
\sum_{k=1}^{n/m_n} \frac{\EE \sup_{\ell=0\ldots m_n}
\left\vert w_{n,i,km_n+\ell}-w_{n,i,km_n}
\right\vert^4}{k^{2-2\eps}}  \le \sum_{k=1}^{n/m_n} \frac{
16 \EE \left\vert w_{n,i,(k+1)m_n}-w_{n,i,km_n}
\right\vert^4}{k^{2-2\eps}}.
$$
Here we used the $L^p$ maximal inequality for martingales,
see  Section 2.2 of \cite{EthierKurtz}.   The expectation is now bounded by a
constant independent of $n,k$, and so is the entire sum, as
required.

\subsection*{The $\beta$-Laguerre case}

Once again \cite{DE1} provides a family of tridiagonal
``$\beta$-Laguerre ensembles", with explicit eigenvalue
densities interpolating between those at $\beta = 1,2,4$.
Take the $n \times n$ bidiagonal random matrix
\begin{equation}
\label{Wmatrix} W_{n, \kappa}^{\beta} =  \frac{1}{\sqrt{\beta}} \left[
\begin{array}{ccccc}
  \widetilde{\chi}_{\beta \kappa} & &&& \\
\chi_{\beta(n-1)} &   \widetilde{\chi}_{\beta(\kappa -1)} &   & &\\
&\ddots & \ddots &  & \\
& & \chi_{\beta 2} & \widetilde{\chi}_{\beta( \kappa - n+2)} &    \\
& & & \chi_{\beta} & \widetilde{\chi}_{\beta( \kappa - n+1)}  \\
 \end{array}\right],
\end{equation}
where  the entries are all independent $\chi$ variables of
the indicated parameter. Here $\kappa \in \RR$ and
necessarily $\kappa
> n-1$. Then, by \cite{DE1}, the eigenvalues of $ (W_{n,
\kappa}^{\beta})^{\dagger}(W_{n, \kappa}^{\beta}) $ have
joint density \eqref{e.wishart.pd}.

While the above puts $\kappa > n-1$,  the obvious duality
reproduces all known real and complex ($\beta=1,2$) results
for any limiting ratio of  dimensions $\kappa \rightarrow
\infty$ and $n \rightarrow \infty$.
This $\beta$ family
generalizes the so-called
``null" Wishart ensembles, distinguishing the important
class of $W \Sigma W^{\dagger}$ type matrices with
non-identity $\Sigma$.  For progress on the spectral edge
of the latter, consult \cite{BaikBen} and
\cite{EK2}.

We now proceed with the proof of Theorem \ref{Lagthm}. It
suffices to prove the claim along a further subsequence of
any given subsequence. This allows us to assume that
$\kappa=\kappa(n)$ is an increasing function of $n$, and
that $n/\kappa(n)\to \vartheta \in[0,1]$. Begin with the
matrix (\ref{Wmatrix}), denoted now simply $W_n$. The
``undressed" ensemble $\beta W_n^\dagger W_n$ has the
processes
\begin{eqnarray}\label{e.xi.diag}
    \widetilde{\chi}_{\b \kappa}^2 + {\chi}_{\b (n-1)}^2, \     \widetilde{\chi}_{\b (\kappa-1)}^2  +  {\chi}_{\b (n-2)}^2,  \
    \widetilde{\chi}_{\b (\kappa-2)}^2 +    {\chi}_{\b (n-3)}^2, \ldots
    \\\label{e.xi.subdiag}
    \widetilde{\chi}_{\b (\kappa-1)} {\chi}_{\b (n-1)},  \
    \widetilde{\chi}_{\b (\kappa-2)} {\chi}_{\b (n-2)},  \
    \widetilde{\chi}_{\b (\kappa-3)} {\chi}_{\b (n-3)}, \ldots
\end{eqnarray}
along the main and off-diagonals, respectively. Up to fist
order, the top left corner of the matrix $W_n^\dagger W_n$
has $n+\kappa$ on the diagonal, and $\sqrt{n\kappa}$
off-diagonal. That is, the top left corner of
$$
\frac{1}{\sqrt{n\kappa}}
 \Bigl( (\sqrt{n}+\sqrt{\kappa})^2 I_n -W_{n}^{\dagger} W_{n}   \Bigr).
$$
is approximately a discrete Laplacian. If time is scaled by
$m_n^{-1}$, then space will have to be scaled by $m_n^2$
for this to converge to the continuum Laplacian. Now the
desired convergence of drift and noise terms determines, up
to constant factors
\begin{equation}
\label{secondm}
 m_{n} = \Bigl( \frac{\sqrt{n \kappa}}{ \sqrt{n} + \sqrt{\kappa}}
 \Bigr)^{2/3},\qquad
H_n= \frac{m_n^2}{\sqrt{n\kappa}}
 \Bigl( (\sqrt{n}+\sqrt{\kappa})^2 I_n -W_{n}^{\dagger}W_{n}   \Bigr).
\end{equation}
Now the $y$'s are defined by formula
\eqref{e.assumptiongrowth}, and are just partial sums of
shifted and scaled versions of \eqref{e.xi.diag} and
\eqref{e.xi.subdiag}.  That is,
\begin{eqnarray*}
 \Da y_{n ,1, k} & = &
\frac{m_n}{\sqrt{n\kappa}} \Bigl( n + \kappa -
 \b^{-1} (\chi_{\b(n-k)}^2  +
 \widetilde{\chi}_{\b(\kappa-k+1)}^2 )
                   \Bigr), \\
       \Da            y_{n ,2,k}  & = &
 \frac{m_n}{\sqrt{n\kappa}} \;2\Bigl(
 \sqrt{n\kappa} -  \b^{-1}  \chi_{\b(n-k)}
\widetilde{\chi}_{\b(\kappa-k)}
 \Bigr).
\end{eqnarray*}
As before, we set $\eta$ to be the expected increments and
$w$ to be the centered $y$.  The $y_{n,i,\cdot},$ $i=1,2$
are again independent increment processes, though they are
not independent of one another. We set
$\gamma=\lim_{n\to\infty}
2\sqrt{n/\kappa}/(\sqrt{n/\kappa}+1)^2\in[0,1/2]$. Then
with $x=m_n/k$ we have
\begin{eqnarray*}
m_n\EE \Da y_{1,n,k} &=&  \gamma\,x+o(1),
\\ m_n\EE (\Da y_{1,n,k})^2&=&\frac{1-\gamma}{2\beta}+o(1),
\qquad m_n\EE (\Da y_{1,n,k}^4)=o(1),
\end{eqnarray*}
uniformly for $k/m_n$ in compacts, so Corollary
\ref{c.either} shows that $y_{n,1}(x)$ converges to the
process $\sqrt{\gamma}/\sqrt{2\beta}\,b_x+\gamma x^2/2$,
whence it is tight. Similarly we get the process convergence
$y_{n,2}(x)\Rightarrow
\sqrt{1-\gamma}/\sqrt{2\beta}\,b_x+(1-\gamma) x^2/2$.

To get the convergence of the sum and \eqref{bounds1}, we
instead consider the process defined by
\begin{eqnarray*}
 \Da  y_{n, k} & = &\frac{m_n}{\sqrt{n\kappa}} \Bigl( (\sqrt{n} + \sqrt{\kappa})^2 -
 \b^{-1} (\chi_{\b(n-k)}  +
 \widetilde{\chi}_{\b(\kappa-k)} )^2,
                   \Bigr)
\end{eqnarray*}
noting that the process $y_{n,k}-y_{n,1,k}-y_{n,2,k}$ is in
sub-scaling and hence converges to the $0$ process in law by a
fourth moment bound. Now $y_{n, k}$  has independent
increments, and the same brand of moment computations already considered
along with Corollary
\ref{c.either} imply that it converges to
$2/\sqrt{\beta}\,b_x+x^2/2$.

Towards tightness, we set $\eta_{n,i,k}=m_n\EE
\Delta{y}_{i,n,k}$. A bit of work shows that for $\beta$
fixed and all $k\ge 1$ $\kappa, n>10$ we have
\begin{eqnarray*}
c_1\frac{k}{m_n} \le    {m_{n}} ( \eta_{n , 1, k}  +
\eta_{n ,2, k} ) \le c_2\frac{k}{m_n}
\end{eqnarray*}
We also have the upper bound $m_n(\eta_{n,2,k})\le 2m_n^2$.
This verifies \eqref{bounds1} and \eqref{bounds2} with
$\bar\eta(x)=x$. The verification of the oscillation bounds
\eqref{bounds3} is identical to the $\beta$-Hermite case.
Indeed, all we used there was that $\sqrt{m_n}w_{n,i,k}$
are martingales whose increments are independent and have
bounded fourth moments.

\bigskip

\noindent{\bf{Acknowledgments }} We would like to thank A.\
Edelman and B.\ Sutton for making earlier versions of
\cite{ES} available to us.  Thanks as well to M.\
Krishnapur, H.P.\  McKean and B.\ Valk\'o for comments, and again
to L.\ Dumaz for help with the proof of \eqref{laure}. 
The work of
B.R.\ was supported in part by NSF grants DMS-0505680
and DMS-0645756, and
that of B.V.\ by a Sloan Foundation fellowship, by the
Canada Research Chair program, and by NSERC and Connaught
research grants.

\sc \bigskip \noindent
Jos{\'e}  A. Ram{\'{\i}}rez\\
Department of Mathematics, \\
 Universidad de Costa Rica,
San Jose  2060, Costa Rica. \\
{\tt  jaramirez@cariari.ucr.ac.cr}

\sc \bigskip \noindent Brian Rider \\  Department of
Mathematics,
\\ University of Colorado at Boulder, Boulder, CO 80309. \\{\tt
brian.rider@colorado.edu}

\sc \bigskip \noindent B\'alint Vir\'ag \\ Departments of
Mathematics and Statistics, \\ University of Toronto, Toronto, ON,
M5S 2E4 Canada.
\\{\tt balint@math.toronto.edu}

\end{document}